\documentclass[a4paper, 11pt]{article}
\usepackage[T1]{fontenc}
\usepackage{times}

\usepackage{amsmath}
\usepackage{amsfonts}
\usepackage{amssymb}
\usepackage{amsthm}

\theoremstyle{plain}
\newtheorem{thm}{Theorem}[section]
\newtheorem{assu}[thm]{Assumption}

\newtheorem{lemma}[thm]{Lemma}
\newtheorem{cor}[thm]{Corollary}
\newtheorem{defi}[thm]{Definition}

\newtheorem{nota}{Notation}

\theoremstyle{remark}
\newtheorem{remark}[thm]{Remark}

\newcommand{\pe}{\mathbb{P}}
\newcommand{\ee}{\mathbb{E}}
\newcommand{\ff}{\mathbb{F}}

\newcommand{\ef}{\mathcal{F}}
\newcommand{\er}{\mathbb{R}}

\newcommand{\en}{\mathcal{N}}

\newcommand{\e}{\mathrm{e}}
\newcommand{\bX}{ {X}}
\newcommand{\bb}{ {b}}
\newcommand{\bW}{ {W}}
\newcommand{\bs}{ {\sigma}}

\newcommand{\bbx}{ {x}}
\newcommand{\bby}{ {y}}
\newcommand{\n}{\noindent}

\newcommand{\I}{{1\kern-0.4em 1}}

\newcommand{\NI}{{|\kern-0.1em |\kern-0.1em |}}
\newcommand{\A}{{\mathcal A}}

\newcommand{\bt}{\begin{thm}}
\newcommand{\et}{\end{thm}}
\renewcommand{\proof}{\n{\bf Proof.\hspace*{2pt}}}

\newcommand{\sr}{\kern-0.3em:\,}

\newcommand{\qedm}{\vspace*{-\baselineskip}\par\hbox to \hsize{\hfill
 \hbox{\vrule width6pt height6pt depth0pt}}
\medskip} %

\begin{document}
\title{Partial Schauder estimates for unbounded solutions of the Kolmogorov equation}

\author{
Andrzej Palczewski\footnote{Faculty of Mathematics, Informatics and Mechanics, University of Warsaw, Banacha 2, 02-097 Warszawa, Poland (e-mail: A.Palczewski@mimuw.edu.pl)}
 }

\date{ }

\maketitle

\begin{abstract}
The goal of the paper is to show, under possibly weak assumptions, that the function given by the Feynman-Kac formula is a classical solution of the associated Kolmogorov equation. We also show that although this solution is unbounded it fulfills Schauder's estimates. The paper collects mostly known results on the regularity of solutions of linear parabolic equations of second order and the regularity of solutions of Ito stochastic differential equations. We aim to reproduce these and related results in a unified fashion.

\bigskip


\end{abstract}
\section{Introduction}

The problem of the existence and regularity of solutions for the Kolmogorov
equation is a particular case of a more general analysis of
second-order parabolic equations. We regard this equation as a
diffusion equation
corresponding to a suitable stochastic differential equation. Our ultimate goal is to show that a functional of a solution of this stochastic equation given by the Feynman-Kac formula solves the Kolmogorov equation. We also obtain Schauder's estimates for polynomially growing solutions of this equation.

Our starting point is $\bX(s)= (X^1(s) , \dots, X^d(s) )$ a $d$-dimensional Ito stochastic process driven by a Wiener process $W(s)$ in a probability space $(\Omega,\ef, \pe )$ with a filtration $\ff= (\ef^t_s)_{t\le s\le T}$ generated by $W(s)$. The notation $\ef^t_s$ means that we are considering the Wiener process $W(s)$ starting at $s= t$ and that the $\sigma$-field $\ef^t_t$ is trivial.
The process $X(s)$ is defined for $s\in [t,T]$ by the system of stochastic differential equations
\begin{equation}\label{eq_111}
d\bX(s) = \bb(s, \bX(s)) ds + \bs(s, \bX(s)) d\bW(s), \quad \bX(t) = \bbx.
\end{equation}
Here $b$ is a $d$-dimensional vector, $\sigma$ a matrix of dimension $d\times
m$, and $\bW(s)$ is an $m$-dimensional Wiener process. To indicate that $\bX(s)$
is a solution of \eqref{eq_111} with initial condition $\bX(t) = \bbx$ we shall write $X^{t,x}(s)$.

Assuming that the coefficients $\bb(s, \bby)$ and
$\bs(s, \bby)$ are deterministic functions we define
the family of elliptic operators
\begin{equation}\label{eq_110}
\A^t u(t,\bbx) = \sum_{i,j =1}^d a_{ij}(t,\bbx)\frac{\partial^2 u}{\partial x_i\partial x_j } + \sum_{i=1}^d b_i(t,\bbx) \frac{\partial u}{\partial x_i},
\end{equation}
where coefficients $a_{ij}$ are given by the formula $a_{ij}\! =\! \frac12\sum_{k=1}^m \sigma_{ik}\sigma_{jk}$. The Kolmogorov equation, which will be the object of our investigation, is then given by the expression
\begin{equation}\label{eq_111a}
D_tu(t,x) + \A^tu(t,x) - c(t,x)u(t,x) + f(t,x) = 0,
\end{equation}
where $x\in\er^d$ and $t\in [0,T]$.

It is well known (cf. Theorem 5.7.6 in \cite{KS1991})
that if equation \eqref{eq_111a} with continuous coefficients and continuous terminal condition possesses a smooth solution of class $C^{1,2}\bigl([0,T]\times\er^d\bigr)$ then
this solution can be represented by the Feynman-Kac formula.

Equally important is the inverse problem: under which conditions the function
given by the Feynman-Kac formula is a solution of equation
\eqref{eq_111a}. A partial answer to this question is given by
Friedman (cf. Theorem 5.6.1 in \cite{F1975}). A more complete analysis can be
found in the book by Krylov (Theorem V.7.4 in \cite{K1995}). In both
mentioned theorems, the authors require a high smoothness of the coefficients
($C^2$ in $x$ with $t$-uniform bounds on derivatives) and the terminal condition
(also $C^2$). On the other hand, the classical theory of parabolic equations in
bounded domains guarantees the existence of smooth solutions of equation
\eqref{eq_111a} for H\"older continuous coefficients and continuous initial
data (cf. Theorem 3.9 in \cite{F1964}). But the picture is more complicated as the stochastic process $X^{t,x}(s)$ that appears in the Feynman-Kac formula is a solution of equation \eqref{eq_111} in which the coefficients can grow linearly with $x$. Hence, the associated equation \eqref{eq_111a} possesses unbounded coefficients, and the results of Theorem 3.9 in \cite{F1964} cannot be straightforwardly applied.

In the last sixty years, parabolic equations with unbounded coefficients in unbounded domains were studied in great detail. After the pioneering works by Aronson and Besala \cite{AB1967}, and Besala \cite{B1975}, the problem was studied by many authors who proved the existence, uniqueness, and regularity of solutions under different hypotheses on the coefficients in spaces of bounded continuous functions but also in $L^p$ spaces on $\er^d$ or on unbounded
domains in $\er^d$ (cf. \cite{A2013}, \cite{AL2011},\cite{AL2016}, \cite{ALP2016}, \cite{BGRT2018}, \cite{CV1987}, \cite{KP2010}, \cite{KLL2010}, \cite{L2011}, \cite{LLS2016}, \cite{L1998}, \cite{M2019}).

An essential difficulty in simultaneous analysis of the stochastic equation \eqref{eq_111} and the parabolic equation \eqref{eq_111a} is caused by the assumed regularity of coefficients. The standard assumption for parabolic problems is H\"older's continuity of coefficients in both variables $(t,x)$ (see the book by Friedman \cite{F1964}). On the other hand, a standard assumption for Ito stochastic equations is Lipschitz continuity in $x$ and continuity (or only measurability) in $t$. To be in accord with the assumptions for stochastic differential equations, we should construct a solution for the parabolic equation \eqref{eq_111a} for coefficients only continuous in
$t$, and H\"older continuous in the space
variables. Then one can expect the Schauder estimates controlling only the H\"older $x$-norms of the derivatives
of solutions. Such estimates for the mentioned weaker assumptions on the coefficients were introduced under the name of partial (or intermediate) Schauder estimates (cf. \cite{B1969}, \cite{K1980}, \cite{L1992}). Then, to obtain the existence and uniqueness of solutions for linear parabolic problems in a bounded domain, it is sufficient to assume the continuity of coefficients in $t $ and H\"older's continuity in $x$, and the continuity of initial and boundary data (cf. \cite{JS2006}). These assumptions are weaker than the standard assumptions for stochastic equations where Lipschitz continuity in $x$ is required, but we are not going to discuss this problem further (for the existence of unique solutions of equation \eqref{eq_111} with local Lipschitz continuity of coefficients see Chapter 1 of the book by Cerrai \cite{C2001}). Successfully, partial Schauder's estimates were extended to bounded solutions of parabolic problems on the whole space (with lower order unbounded terms in \cite{KP2010} and all unbounded coefficients in \cite{L2011}). Partial Schauder's estimates for unbounded solutions (in weighted H\"older spaces) were obtained in \cite{L2000} but only for bounded coefficients. Our goal in this paper is to obtain these estimates in the whole space for unbounded solutions and unbounded coefficients. These results are an extension of results in \cite{L2000} and \cite{L2011}.

Time regularity of the coefficients in the Kolmogorov equation requires an additional comment. There are papers assuming only $(t,x)$ measurability of the coefficients and H\"older's regularity in the spatial variables (cf. \cite{BL2005}, \cite{L2011}). Under these assumptions it can be proved that the Kolmogorov equation is fulfilled in a generalized sense, i.e., the solution $u(t,x)$ is a Lipschitz continuous function in $t$ and fulfills the equation outside a negligible set. We are not going to explore this generalization, as our goal is the construction of classical solutions of equation \eqref{eq_111a}, which requires the assumption of $t$-continuity of the coefficients.
Then under suitable assumptions on the space regularity, we get strong solutions of equation \eqref{eq_111} and $C^{1,2}\bigl((0,T)\times\er^d\bigr)$ solutions of equation \eqref{eq_111a}.

The rest of the paper is organized as follows. In Section 2 we prove under relatively weak assumptions that the function given by the Feynman-Kac formula is a classical solution of the Kolmogorov equation (Theorem \ref{th_142}). We also obtain estimates of the spatial derivatives of this solution (Theorem \ref{th_144} and Corollary \ref{cor_141}). In Section 3 we extend the estimates of Section 2 to partial Schauder's estimates in   spaces of polynomially growing functions (Theorem \ref{th_166}).

\section{Equations of Kolmogorov type}

We will now concentrate on the probabilistic representation of the solutions to the Cauchy problem of Kolmogorov's equation
\begin{equation}\label{eq_142}
\begin{split}
&D_tu(t,x) + \A^tu(t,x) - c(t,x)u(t,x) + f(t,x) = 0, \ (t,x) \in [0,T)\times\er^d,\\
&u(T,x) = h(x), \quad x\in \er^d,
\end{split}
\end{equation}
where $\A^t$ given by equation \eqref{eq_110} is associated with the stochastic differential equation
\begin{equation}\label{eq_141}
d\bX(s) = \bb(s, \bX(s)) ds + \bs(s, \bX(s)) d\bW(s),\quad \bX(t) = \bbx.
\end{equation}
About the coefficients of equation \eqref{eq_141} we make the following assumptions.

\begin{assu}\label{as_121}
Let the coefficients in equation \eqref{eq_141} be deterministic, continuous functions in $(s,\bby)$ and such that:
\begin{enumerate}
\item[(A1)] $|\bb(s,\bby) - \bb(s,\bby')| + |\bs(s,\bby) - \bs(s,\bby')| \le C |\bby-\bby'|$;
\item[(A2)] $|\bb(s, \bby)| +|\bs(s, \bby)| \le C (1+|\bby|)$;
\end{enumerate}
for $\bby, \bby' \in \er^d$ and $s \in [t,T]$, where $|\cdot|$ denotes the Euclidean norm in $\er^d$ and $C$ is a finite, positive constant depending on $T$.

\n
We also assume that $x$ is a deterministic vector in $\er^d$.
\end{assu}

Under the above assumptions, equation \eqref{eq_141} possesses a unique strong solution $\bX^{t,x}(s)$.

Our goal is to prove that a solution of the Cauchy problem \eqref{eq_142} is given by the Feynman-Kac formula
\begin{equation}\label{eq_143}
\begin{split}
u(t,x) = \, &\ee\biggl(h\bigl(X^{t,x}(T)\bigr)\exp\Bigl(-\int_t^T c\bigl(r, X^{t,x}(r)\bigr)dr\Bigr)\biggr.\\
& \biggl.+ \int_t^T f\bigl(s,X^{t,x}(s)\bigr)\exp\Bigl(-\int_t^s c\bigl(r, X^{t,x}(r)\bigr)dr\Bigr)ds\biggr),
\end{split}
\end{equation}
where $\bX^{t,x}(s)$ is a solution of equation \eqref{eq_141}.

In the rest of the paper, we will use the following notation.

\begin{nota}%
\n
\begin{enumerate}
\item The space of functions $f\sr \mathcal{B}\to \er$ which are bounded and uniformly continuous on $\mathcal{B}$ will be denoted $BC(\mathcal{B})$. This is a Banach space with the norm
\[
\|f\|_{BC(\mathcal{B})} = \sup_{z\in \mathcal{B}} |f(z)|.
\]
Here $\mathcal{B}$ can be a subset of $[0,T]\times \er^d$ or a subset of $\er^d$.
\item The space $BC^{p}(U)$, $p\ge 1$, where $U\subset \er^d$, is a space of functions $f\sr U\to \er$ with continuous, bounded derivatives $D^\alpha f(x)\in BC(U)$, for $1\le |\alpha| \le p$. This space is a Banach space with the norm
\[
\|f\|_{BC^{p}(U)} = \sum_{ |\alpha|\le p }\sup_{x\in U} |D^\alpha f(x)|.
\]

\item Let $U$ an open, bounded domain in $\er^d$. We will always assume that an open bounded set $U$ is compactly embedded in $B(R)$, an open ball in $\er^d$ of radius $R$, i.e., $\bar U \subset B(R)$. For an open, bounded set $U$ we define the H\"older space $H^\beta(U)$, a space of functions $f\sr U\to \er$ which are in $BC(U)$ and are $\beta$-H\"older continuous, that is
\[
[f]_\beta = \sup_{ \genfrac{}{}{0pt}{}{x,x' \in U}{ x\neq x'} }
\biggl(\frac{| f(x) - f(x')|}{|x - x'|^\beta}\biggl) < + \infty.
\]
The H\"older space $H^{p+\beta}(U)$ is a space of functions $f(x)$ with continuous derivatives $D^\alpha f(x)\in BC(U)$, for $0\le |\alpha| \le p$, and $\beta$-H\"older continuous derivatives $D^\alpha f(x)$ for $|\alpha| = p$.

The space $H^{p+\beta}(U)$ is a Banach space for all $p\ge 0$ and $\beta\in (0,1)$ with the norm
\[
\|f\|_{H^{p+ \beta}(U)} = \sum_{ |\alpha|\le p }\sup_{x\in U} |D^\alpha f(x)| + \sum_{ |\alpha|= p } \bigl[D^{\alpha} f(x)\bigr]_\beta.
\]
Analogously, we define the space $H^{p+\beta}(\er^d)$.
\item The H\"older space $H^{p+\beta}_{\mathit{loc}}(\er^d)$ is a space of functions $f\sr \er^d \to \er$, which for any open, bounded set $U\subset \er^d$ belong to $H^{p+\beta}(U)$.
\item Let $P(x) = 1+|x|^{2q}$ denote for an integer $q \ge 0$ the polynomial weight in $\er^d$. The space of functions $f\sr \er^d\to \er$ such that $ f P^{-1}\in BC^p(\er^d)$, $p\ge 0$, will be denoted $BC^p_P(\er^d)$. This is a Banach space with the norm
\[
\|f\|_{BC^p_P(\er^d)} = \sum_{ |\alpha|\le p } \sup_{x\in \er^d} \bigg|D^\alpha\frac{f(x)}{P(x)}\bigg|.
\]
\item $H^{p+\beta}_P(\er^d)$ denotes the space of functions $f\sr \er^d\to \er$ such that $f P^{-1}\in H^{p+\beta}(\er^d)$. This is a Banach space with the norm
\[
\|f\|_{H^{p+ \beta}_P(\er^d)} = \sum_{ |\alpha|\le p }\sup_{x\in\er^d} \bigg|D^\alpha \frac{f(x)}{P(x)}\bigg| + \sum_{ |\alpha|= p } \biggl[D^{\alpha}\frac{f(x)}{P(x)}\biggr]_\beta.
\]
\end{enumerate}
\end{nota}

By simple computations, we obtain the following lemma that gives an equivalent definition of the norm in $H^{p+\beta}_P(\er^d)$ (this is Lemma 2.5 of \cite{L2000}).
\begin{lemma}
For $f\in H^{p+\beta}_P(\er^d)$ the norm $\|f\|_{H^{p+ \beta}_P(\er^d)}$ is equivalent to the norm
\[
\NI f\NI_{H^{p+ \beta}_P(\er^d)} = \sum_{ |\alpha|\le p }\sup_{x\in\er^d} \bigg|\frac{D^\alpha f(x)}{P(x)}\bigg| + \sum_{ |\alpha|= p } \biggl[\frac{D^{\alpha}f(x)}{P(x)}\biggr]_\beta.
\]
\end{lemma}

We begin with the proof that under suitable assumptions, function $u(t,x)$ given by equation \eqref{eq_143} is a continuous function for $(t,x) \in [0,T]\times \er^d$. Results of this type are proved in many papers on the Kolmogorov equation with unbounded coefficients, but under very restrictive assumptions on the regularity of the coefficients (cf. \cite{AL2015}, \cite{KLL2010}, \cite{L2011}). In the paper by Rubio \cite{Rubio2}, the continuity of $u(t,x)$ was proved assuming $t$-H\"older and $x$-Lipschitz continuity of the coefficients and some regularity of initial data. The following theorem is much stronger as we assume only the continuity of coefficients and  initial data, and H\"older's continuity in the spatial variables of the coefficients.

\bt\label{th_141} %
Let the coefficients of equation \eqref{eq_141} fulfill the conditions of Assumptions \ref{as_121}, i.e., $\bb(t,\bbx)$ and $\bs(t,\bbx)$ be $(t,x)$-continuous for $(t,x) \in [0,T]\times\er^d$, and with respect to $x$ be Lipschitz continuous with linear growth.

The coefficients of equation \eqref{eq_142} fulfill the following assumptions:
\begin{enumerate}
\item $c\sr[0,T] \times \er^d\to \er$\, and $f\sr[0,T] \times \er^d\to \er$\, are $(t,x)$-continuous and locally H\"older continuous with respect to $x$, i.e., for an open bounded set $U$ such that $\bar U \subset B(R)$, and $t\in(0,T)$, $c(t, \cdot)$, and $f(t,\cdot)$ belong to $H^{\beta} (U )$, $0< \beta<1$, and
\[
\|c(t, \cdot)\|_{H^{\beta}(U)} + \|f(t, \cdot)\|_{H^{\beta}(U)} \le C(T,R).
\]
\item $h\sr \er^d \to \er$\, is $x$-continuous.
\item There is $c_0 >0$ such that $c(t,x) \ge c_0$ for $(t,x)\in [0,T]\times\er^d$.
\item There is an integer exponent $q \ge 1$ such that:
\begin{enumerate}
\item[ ] $|f(t,x)| \le C(T) \bigl(1+|x|^{2q}\bigr)$ for$(t,x)\in [0,T]\times\er^d$.
\item[ ] $|h( x)| \le C \bigl(1+|x|^{2q} \bigr)$ for $x\in \er^d$.
\end{enumerate}

\end{enumerate}
Under the above assumptions the function $u(t,x)$ given by equation \eqref{eq_143} belongs to $C\bigl([0,T]\times\er^d\bigr)$. Furthermore
\[
\sup_{0\le t\le T} |u(t,x)| \le C(T) \bigl(1+ |x|^{2q}\bigr), \quad x\in\er^d,
\]
where $q$ is the exponent   defining the rate of growth of $f(t,x)$ and $h(x)$.
\et
\proof
The proof explores many ideas of the proof of Theorem 3.1 of \cite{Rubio2}. We are repeating this proof not only for the reader's convenience but to point out the differences that make our proof possible under much weaker assumptions.

The proof of continuity of $u(t,x)$ on the interval $[0,T]$ requires special treatment of boundary points. To avoid these problems, we consider the interval $(-\gamma, T]$ for $\gamma$ positive but small. This overcomes the problem of the boundary point $t=0$. The point $t=T$ still requires special consideration. Let $(t,x) \in (-\gamma, T]\times \er^d$ and $\{(t_n,x_n)\}_{n\in \en}$ be a sequence of points $(t_n,x_n) \in(-\gamma, T]\times \er^d$ converging to $(t,x)$. We assume that $t_n \nearrow t$, then the point $t=T$ can be treated in the same way as interior points.
This approach simplifies computations and puts no restriction on the generality of considerations.

The proof of convergence
\[
u(t_n,x_n) \to u(t,x)
\]
can be reduced by formula \eqref{eq_143} to the following
\begin{align*}
\ee\biggl(&h\bigl(X^{t_n,x_n}(T)\bigr)\exp\Bigl(-\int_{t_n}^T c\bigl(r, X^{t_n,x_n}(r)\bigr)dr\Bigr)\biggr.\\
&\biggl. - h\bigl(X^{t,x}(T)\bigr)\exp\Bigl(-\int_t^T c\bigl(r, X^{t,x}(r)\bigr)dr\Bigr)\biggr)\\
+ \ee \biggl(& \int_{t_n}^T f\bigl(s,X^{t_n,x_n}(s)\bigr)\exp\Bigl(-\int_{t_n}^s c\bigl(r, X^{t_n,x_n}(r)\bigr)dr\Bigr)ds\biggr.\\
&\biggl.- \int_t^T f\bigl(s,X^{t,x}(s)\bigr)\exp\Bigl(-\int_t^s c\bigl(r, X^{t,x}(r)\bigr)dr\Bigr)ds\biggr) \to 0.
\end{align*}
The strategy of the proof is to show that the expectations of random variables
\begin{align*}
v^n_1 =\, &\biggl|h\bigl(X^{t_n,x_n}(T)\bigr)\exp\Bigl(-\int_{t_n}^T c\bigl(r, X^{t_n,x_n}(r)\bigr)dr\Bigr)\biggr.\\
&\biggl. - h\bigl(X^{t,x}(T)\bigr)\exp\Bigl(-\int_t^T c\bigl(r, X^{t,x}(r)\bigr)dr\Bigr)\biggr|
\end{align*}
and
\begin{align*}
v^n_2 = \, & \biggl| \int_{t_n}^T f\bigl(s,X^{t_n,x_n}(s)\bigr)\exp\Bigl(-\int_{t_n}^s c\bigl(r, X^{t_n,x_n}(r)\bigr)dr\Bigr)ds\biggr.\\
&\biggl.- \int_t^T f\bigl(s,X^{t,x}(s)\bigr)\exp\Bigl(-\int_t^s c\bigl(r, X^{t,x}(r)\bigr)dr\Bigr)ds\biggr|
\end{align*}
converge to zero.

An essential step in the proof is to show that variables $v^n_1$ and $v^n_2$ are uniformly integrable. We obtain the required uniform integrability by proving the square integrability of each of these variables. 

For $v^n_1$ we have
\begin{align*}
\int_\Omega |v^n_1|^2 d\pe &\le 2 \int_\Omega \biggl(|h\bigl(X^{t_n,x_n}(T)\bigr)| \exp\Bigl(-\int_{t_n}^T c\bigl(r, X^{t_n,x_n}(r)\bigr)dr\Bigr)\biggr)^2 d\pe\\
&+ 2 \int_\Omega \biggl(|h\bigl(X^{t,x}(T)\bigr)| \exp\Bigl(-\int_{t}^T c\bigl(r, X^{t,x}(r)\bigr)dr\Bigr)\biggr)^2 d\pe\\
\le & C \int_\Omega \Bigl(\bigl(1 + |X^{t,x}(T)|^{4q}\bigr) + \bigl(1 + |X^{t_n,x_n}(T)|^{4q}\bigr)\Bigr) d\pe \le C (1+|x|^{4q}).
\end{align*}
Thus $v^n_1$ is uniformly integrable.

To prove that $\ee(v^n_1)\to 0$ as $n\to \infty$ we define the set
\begin{align*}
B_{M,\eta} =&
\{\omega\in\Omega\sr \sup_{s\in[t,T]} |X^{t,x}(\omega, s)|\le M, \sup_{s\in[t_n,T]} |X^{t_n,x_n}(\omega, s)|\le M, \\
&\sup_{s\in[t_n,T]} |X^{t,x}(\omega, s) - X^{t_n,x_n}(\omega, s)|\le \eta \}.
\end{align*}
To avoid problem with the comparison of processes $X^{t_n,x_n}$ and $X^{t,x}$ on the interval $(t_n, t)$ we extend $X^{t,x}$ on this interval putting $X^{t,x}(s) = x$ for $s\in[t_n,t)$.

It follows from the uniform integrability of $v^n_1$ that taking $M$ sufficiently large and $\eta$ sufficiently small we obtain
\[
\int_{\Omega\setminus B_{M,\eta}} v^n_1 d\pe \le \delta,
\]
where $\delta$ is a small number. To estimate the integral over $B_{M,\eta}$, we split $v^n_1$
\begin{equation}\label{eq_144}
\begin{split}
v^n_1 \le &\, \Big|h\bigl(X^{t_n,x_n}(T)\bigr) - h\bigl(X^{t,x}(T)\bigr)\Big| \exp\Bigl(-\int_{t_n}^T c\bigl(r, X^{t_n,x_n}(r)\bigr)dr\Bigr)\\
& + \Big|h\bigl(X^{t,x}(T)\bigr)\Big|\\
&\times\bigg|\exp\Bigl(-\int_{t_n}^T c\bigl(r, X^{t_n,x_n}(r)\bigr)dr\Bigr) - \exp\Bigl(-\int_t^T c\bigl(r, X^{t,x}(r)\bigr)dr\Bigr)\bigg|.
\end{split}
\end{equation}

For the mean value of the first term in the right hand side of \eqref{eq_144}, we have
\begin{align*}
\ee\biggl(& \Big|h\bigl(X^{t_n,x_n}(T)\bigr) - h\bigl(X^{t,x}(T)\bigr)\Big|\exp\Bigl(-\int_t^T c\bigl(r, X^{t,x}(r)\bigr)dr\Bigr)\biggr)\\
&\le \ee \Big|h\bigl(X^{t_n,x_n}(T)\bigr) - h\bigl(X^{t,x}(T)\bigr)\Big| \\
&\le C \Bigl(1 + \ee\bigl|X^{t_n,x_n}(T)\bigr|^{2q} + \ee\bigl|X^{t,x}(T)\bigr|^{2q} \Bigr)
\le C\bigl(1 +|x|^{2q}\bigr).
\end{align*}
Since $X^{t_n,x_n}(T)$ converges in probability to $X^{t,x}(T)$ (cf. Theorem II.2.1 in \cite{K1984} and Lemma III.6.13 in \cite {K1995}) then by Lemma III.6.13 in \cite {K1995} we have
\[
\ee\biggl( \Big|h\bigl(X^{t_n,x_n}(T)\bigr) - h\bigl(X^{t,x}(T)\bigr)\Big|\exp\Bigl(-\int_t^T c\bigl(r, X^{t,x}(r)\bigr)dr\Bigr)\biggr)\to 0.
\]

Choosing $n$ large enough, we obtain from the above
\[
\int_{B_{M,\eta}}\biggl( \Big|h\bigl(X^{t_n,x_n}(T)\bigr) - h\bigl(X^{t,x}(T)\bigr)\Big|\exp\Bigl(-\int_t^T c\bigl(r, X^{t,x}(r)\bigr)dr\Bigr)\biggr)d\pe \le \delta/2.
\]

The estimation of the integral of the second term in the right hand side of \eqref{eq_144} we begin with the estimation of the difference of exponents
\begin{align*}
&\bigg|\exp\Bigl(-\int_{t_n}^T c\bigl(r, X^{t_n,x_n}(r)\bigr)dr\Bigr) - \exp\Bigl(-\int_t^T c\bigl(r, X^{t,x}(r)\bigr)dr\Bigr)\bigg|\\
&\le \bigg|\exp\Bigl(-\int_{t_n}^T\Bigl( c\bigl(r, X^{t_n,x_n}(r)\bigr) - c\bigl(r, X^{t,x}(r)\bigr)\Bigr)dr\Bigr) -1 \bigg|\\
&\le \biggl( \exp\Bigl(\int_{t_n}^T\Big| c\bigl(r, X^{t_n,x_n}(r)\bigr) - c\bigl(r, X^{t,x}(r)\bigr)\Big|dr\Bigr) -1 \biggr)\\
&\le \biggl( \exp\Bigl(\int_{t_n}^T C\big| X^{t_n,x_n}(r) - X^{t,x}(r)\big|^\beta dr\Bigr) -1 \biggr)\\
&\le \e^{CT \eta^\beta} -1 \le C \theta T \eta^\beta,
\end{align*}
where we used H\"older's continuity of $c(t,x)$ and the mean value theorem which introduced $\theta\in [0,1]$.
Then for a sufficiently small $\eta$, we obtain
\begin{align*}
&\int_{B_{M,\eta}} \biggl(\big|h\bigl(X^{t,x}(T)\bigr)\big|\biggr.\\
\biggl. &\times \bigg|\exp\Bigl(-\int_{t_n}^T c\bigl(r, X^{t_n,x_n}(r)\bigr)dr\Bigr) - \exp\Bigl(-\int_t^T c\bigl(r, X^{t,x}(r)\bigr)dr\Bigr)\bigg|\biggr) d\pe \le \delta/2.
\end{align*}

Collecting the above estimates, we have
\[
\ee(v^n_1) = \int_\Omega v^n_1 d\pe \le 2 \delta.
\]

Similarly to $v^n_1$, we prove the uniform integrability of $v^n_2$. To this end, we compute
\begin{align*}
\int_\Omega |v^n_2|^2 &d\pe \le 2 \int_\Omega \biggl(\int_{t_n}^T f\bigl(s, X^{t_n,x_n}(s)\bigr) \exp\Bigl(-\int_{t_n}^s c\bigl(r, X^{t_n,x_n}(r)\bigr)dr\Bigr)ds\biggr)^2d\pe \\
&+ 2 \int_\Omega \biggl(\int_{t}^T f\bigl(s, X^{t,x}(s)\bigr) \exp\Bigl(-\int_{t}^s c\bigl(r, X^{t,x}(r)\bigr)dr\Bigr)ds\biggr)^2d\pe \\
\le & C \int_\Omega \int_{t_n}^T \Bigl(\bigl(1 +\sup_{s\in[t,T]} |X^{t,x}(s)|^{4q}\bigr) + \bigl(1 + \sup_{s\in[t_n,T]} |X^{t_n,x_n}(s)|^{4q}\bigr)\Bigr)ds d\pe\\
\le & C T \int_\Omega \Bigl(2 +\sup_{s\in[t,T]} |X^{t,x}(s)|^{4q} + \sup_{s\in[t_n,T]} |X^{t_n,x_n}(s)|^{4q} \Bigr) d\pe\\
\le & C T \bigl(2 + |x|^{4q} + |x_n|^{4q}\bigr) \le CT (1+|x|^{4q}).
\end{align*}
We define the set $B_{M,\eta}$ in a similar way as for $v^n_1$ and select $M$ and $\eta$ such that
\[
\int_{\Omega\setminus B_{M,\eta}} v^n_2 d\pe \le \delta.
\]
To estimate the integral over $B_{M,\eta}$ we split $v^n_2$ as follows
\begin{equation}\label{eq_145}
\begin{split}
v^n_2 \le &\, \int_{t_n}^T \Big|f\bigl(s, X^{t_n,x_n}(s)\bigr) - f\bigl(s, X^{t,x}(s)\bigr)\Big| \exp\Bigl(-\int_{t_n}^s c\bigl(r, X^{t_n,x_n}(r)\bigr)dr\Bigr)ds\\
& + \int_{t_n}^T\Big|f\bigl(s,X^{t,x}(T)\bigr)\Big|\\
&\times\bigg|\exp\Bigl(-\int_{t_n}^s c\bigl(r, X^{t_n,x_n}(r)\bigr)dr\Bigr) - \exp\Bigl(-\int_t^s c\bigl(r, X^{t,x}(r)\bigr)dr\Bigr)\bigg|ds
\end{split}
\end{equation}
For the first integral in the right hand side of \eqref{eq_145}, we apply H\"older's continuity of $f(t,x)$ to obtain
\begin{align*}
\int_{B_{M,\eta}} &\int_{t_n}^T \!\!\Big|f\bigl(s, X^{t_n,x_n}(s)\bigr) - f\bigl(s, X^{t,x}(s)\bigr)\Big| \exp\Bigl(\!-\!\int_{t_n}^s c\bigl(r, X^{t_n,x_n}(r)\bigr)dr\Bigr)ds d\pe\\
&\le \, C \int_{B_{M,\eta}}\int_{t_n}^T \Big|f\bigl(s, X^{t_n,x_n}(s)\bigr) - f\bigl(s, X^{t,x}(s)\bigr)\Big| ds d\pe\\
&\le \, C \int_{B_{M,\eta}}\int_{t_n}^T \big| X^{t_n,x_n}(s)- X^{t ,x }(s)\big|^\beta dsd\pe \le \, C T\eta^\beta.
\end{align*}

To estimate the second integral in the right hand side of \eqref{eq_145} for $t_n \le s\le T$, we evaluate,   like for $v^n_1$,  the difference of exponents
\begin{align*}
\bigg|\exp\Bigl(-\int_{t_n}^s c\bigl(r, X^{t_n,x_n}(r)\bigr)dr\Bigr) - \exp\Bigl(-\int_t^s c\bigl(r, X^{t,x}(r)\bigr)dr\Bigr)\bigg|
\le C \theta s \eta^\beta,
\end{align*}

Then, for a sufficiently small $\eta$, we obtain for the second integral in the right hand side of \eqref{eq_145}
\begin{align*}
\int_{B_{M,\eta}} & \int_{t_n}^T\Big|f\bigl(s,X^{t,x}(T)\bigr)\Big| \\
&\times\bigg|\exp\Bigl(-\int_{t_n}^s c\bigl(r, X^{t_n,x_n}(r)\bigr)dr\Bigr) - \exp\Bigl(-\int_t^s c\bigl(r, X^{t,x}(r)\bigr)dr\Bigr)\bigg|dsd\pe \\
&\le\, C \eta^\beta \,\int_{B_{M,\eta}} \int_{t_n}^T s \Big|f\bigl(s,X^{t,x}(s)\bigr)\Big|ds \\
&\le\, C \eta^\beta \,\int_{B_{M,\eta}}\int_{t_n}^T s \Bigl(1 + \sup_{s\in [t,T]}\big|X^{t,x}(s)\big|^{2q} \Bigr)dsd\pe
\le\, C T^2 \eta^\beta \bigl(1 +|x|^{2q} \bigr).
\end{align*}
By the above estimates, we have for a sufficiently small $\eta$
\begin{align*}
\ee(v^n_2) = \int_\Omega v^n_2 d\pe \le \int_{B_{M,\eta}} v^n_2 d\pe + \delta \le C \eta^\beta(1+|x|^{2q}) +\delta \le 2 \delta.
\end{align*}
This proves that $u(t,x)$ is continuous in $ (-\gamma, T]\times \er^d$, hence also in $ [0, T]\times \er^d$.

The growth estimate of $u(t,x)$ can be easily obtained from equation \eqref{eq_143}. For the first term in \eqref{eq_143}, we have
\begin{align*}
\bigg| \ee\biggl(& h\bigl(X^{t,x}(T)\bigr) \exp\Bigl(-\int_t^T c\bigl(r, X^{t,x}(r)\bigr)dr\Bigr)\biggr)\bigg|\\
&\le \ee \Big|h\bigl(X^{t,x}(T)\bigr)\Big|
\le C \Bigl(1 + \ee\bigl|X^{t,x}(T)\bigr|^{2q} \Bigr)
\le C\bigl(1 +|x|^{2q}\bigr),
\end{align*}
where the constant $C$ depends on $T$, $q$, and the constant in the growth estimate of $h(x)$.

For the second term in \eqref{eq_143} we obtain the estimate
\begin{align*}
\bigg|&\ee\biggl(\int_{t}^T f\bigl(s, X^{t,x}(s)\bigr) \exp\Bigl(-\int_{t}^s c\bigl(r, X^{t,x}(r)\bigr)dr\Bigr)ds\biggr)\bigg|\\
&\le \, C \int_{t}^T \ee\Big|f\bigl(s, X^{t,x}(s)\bigr)\Big| ds
\le \, C \int_{t }^T \Bigl(1 + \ee\big|X^{t ,x }(s)\big|^{2q}\Bigr)ds \le \, C \bigl(1 + |x|^{2q}\bigr),
\end{align*}
where the constant $C$ depends on $T$, $q$, and the constant in the growth estimate of $f(t,x)$.

The above estimates give
\[
\sup_{0\le t\le T} |u(t,x)| \le \, C\bigl(1+|x|^{2q}\bigr),
\]
where the constant $C$ depends on $T$, $q$, and the constants in the growth estimates of $f(t,x)$ and $h(x)$.
\qed
\medskip

The continuous function $u(t,x)$ defined by equation \eqref{eq_143} is not a solution of the Kolmogorov equation \eqref{eq_142}. But under slightly stronger assumptions on the coefficients of equation \eqref{eq_142}, this function is already in $C^{1,2}\bigl((0,T)\times\er^d\bigr)$ and is a solution of the Cauchy problem \eqref{eq_142}.

\bt\label{th_142} %
Let the conditions of Assumption \ref{as_121} be fulfilled and operator $\A^t$ be uniformly elliptic, i.e., there exists $\delta(T) >0$ such that for $(t,x) \in [0,T]\times\er^d$ and $\xi\in \er^d\setminus \{0\}$
\[
\sum_{i,j =1}^d a_{ij}(t,x) \xi_i\xi_j \ge \delta(T) |\xi|^2.
\]
Let $h\sr \er^d \to \er$\, belong to $BC(\er^d)$, $c\sr[0,T] \times \er^d\to \er$\, and $f\sr[0,T] \times \er^d\to \er$\, be $(t,x)$-continuous and locally H\"older continuous with respect to $x$, i.e., for an open, bounded set $U$ such that $\bar U \subset B(R)$, and $t\in(0,T)$, $c(t, \cdot)$, and $f(t,\cdot)$ belong to $H^{\beta} (U )$, $0< \beta<1$, and
\[
\|c(t, \cdot)\|_{H^{\beta}(U)} + \|f(t, \cdot)\|_{H^{\beta}(U)} \le C(T,R).
\]
In addition, let $c(t,x) \ge c_0>0$\, for $(t,x)\in [0,T]\times\er^d$.

Then $u(t,x)$ defined by equation \eqref{eq_143} is such that $u(t,\cdot) \in H^{2+ \beta} (U)$ for $t\in (0,T)$.
Furthermore, $u(t,x)$ belongs to $C^{1,2}\bigl((0,T)\times\er^d\bigr)$ and is a unique classical solution of the Cauchy problem \eqref{eq_142}.
\et

\proof
For $(t,x) \in(t_1,t_2)\times U$, where $(t_1,t_2) \subset [0, T]$ and $U$ is an open set with a regular boundary $\partial U$ of class $C^2$, consider the following Cauchy problem with the coefficients fulfilling the assumptions of the theorem
\begin{equation}\label{eq_146}
\begin{split}
&D_tv(t,x) + \A^tv(t,x) - c(t,x)v(t,x) + f(t,x) = 0, \ (t,x) \in (t_1,t_2)\times U,\\
&v(t_2,x) = u(t_2,x), \quad x\in U,\\
&v(t,x) = u(t,x), \quad (t,x)\in (t_1,t_2) \times \partial U,
\end{split}
\end{equation}
where $u(t,x)$ defined by equation \eqref{eq_143} is a continuous function by Theorem \ref{th_141}.

To obtain a regular solution of problem \eqref{eq_146}, let us observe that due to conditions (A1) and (A2) of Assumption \ref{as_121} the coefficients $a_{ij}$ and $b_i$ as functions of $x$ belong to $H^{\beta} (U )$ and
\[
\|a_{ij}(t, \cdot)\|_{H^{\beta}(U)} + \|b_i(t, \cdot)\|_{H^{\beta}(U)} \le C(T,R).
\]
Then,
we approximate $u(t,x)$ by a sequence of continuous functions $u_n(t,x)$ which are H\"older continuous in $x$ with exponent $\beta$. Under these conditions, we can apply Theorem 16.1 of \cite{L1992}. According to this theorem, there exist unique solutions $v_n(t,x)$ of problem \eqref{eq_146} with initial/boundary data $u_n(t,x)$ such that $v_n(t,x) \in C^{1,2} \bigl((t_1, t_2)\times U\bigr)$. Moreover, $v_n(t, \cdot) \in H^{2+ \beta} ( U )$ for all $t\in (t_1,t_2)$. Similarly like in the proof of Lemma 2.6 in \cite{JS2006}, we can pass to the limit with solutions $v_n(t,x)$ using the maximum principle and the condition $c(t,x) \ge c_0 >0$, and obtain a unique solution $v(t,x) \in C^{1,2} \bigl((t_1, t_2)\times U\bigr)$ of problem \eqref{eq_146} with continuous initial/boundary data $u(t,x)$.

To obtain a stochastic representation of solution $v(t,x)$ under the weak assumption of the theorem, we apply Theorem 7.4.2 of \cite{YZ1999} (see also Theorem IV.4.5 of \cite{G2016}). Then we get
\begin{equation}\label{eq_147}
\begin{split}
v(t,x) = \, &\ee_{t,x}\biggl(u\bigl((\tau,X^{t,x}(\tau)\bigr)\exp\Bigl(-\int_t^\tau c\bigl(r, X^{t,x}(r)\bigr)dr\Bigr)\biggr)\\
& + \ee_{t,x}\biggl( \int_t^\tau f\bigl(s,X^{t,x}(s)\bigr)\exp\Bigl(-\int_t^s c\bigl(r, X^{t,x}(r)\bigr)dr\Bigr)ds\biggr),
\end{split}
\end{equation}
where $X^{t,x}(s)$ is a solution of the stochastic equation \eqref{eq_141} with initial condition $X(t) = x$ and $\tau = T\wedge\tau_\varepsilon$ with $\tau_\varepsilon$ the first exit time from $U$, i.e., $\tau_\varepsilon = \inf\{s>0\sr X^{t,x}(s) \notin \bar U\}$ ($\ee(\tau_\varepsilon)<\infty$ due to the proof of Theorem 5.5.1 in \cite{F1975}).

Let $w(t,x)$ denote the restriction of $u(t,x)$ defined by equation \eqref{eq_143} to domain $U$. It is clear that $w\bigl(s, X^{t,x}(s)\bigr) = u\bigl(s, X^{t,x}(s)\bigr)$ for $0\le s\le \tau$ and $X(t) = x$ with $x\in U$. Hence
\[
w\bigl(s, X^{t,x}(s)\bigr) = \ee\bigl(u\bigl(s, X^{t,x}(s)\bigr)\big|\ef^t_\tau\bigr)
\]
and
\begin{equation}\label{eq_148}
\begin{split}
w&(t,x) = \, \ee_{t,x}\Biggl(\ee\biggl(h\bigl(X^{t,x}(T)\bigr)\exp\Bigl(-\int_t^T c\bigl(r, X^{t,x}(r)\bigr)dr\Bigr)\bigg|\ef^t_\tau\biggr) \Biggr)\\
& + \ee_{t,x}\Biggl(\ee\biggl(\int_t^T f\bigl(s,X^{t,x}(s)\bigr)\exp\Bigl(-\int_t^s c\bigl(r, X^{t,x}(r)\bigr)dr\Bigr)ds\bigg|\ef^t_\tau\biggr)\Biggr).
\end{split}
\end{equation}

Using the Markov property of process $X(s)$, we obtain for the first integral
\begin{align*}
\ee_{t,x}\Biggl(&\ee\biggl(h\bigl(X^{t,x}(T)\bigr)\exp\Bigl(-\int_t^T c\bigl(r, X^{t,x}(r)\bigr)dr\Bigr)\bigg|\ef^t_\tau\biggr) \Biggr)\\
=\, &\ee_{t,x}\Biggl(\exp\Bigl(-\int_t^\tau c\bigl(r, X^{t,x}(r)\bigr)dr\Bigr)\biggr.\\
\biggl.&\times\ee\biggl(h\bigl(X^{t,x}(T)\bigr)\exp\Bigl(-\int_\tau^T c\bigl(r, X^{t,x}(r)\bigr)dr\Bigr)\bigg|\ef^t_\tau\biggr) \Biggr)\\
=\, &\ee_{t,x}\Biggl(\exp\Bigl(-\int_t^\tau c\bigl(r, X^{t,x}(r)\bigr)dr\Bigr)\biggr.\\
\biggl.&\times\ee_{\tau,X(\tau)}\biggl(h\bigl(X^{\tau,X(\tau)}(T)\bigr)\exp\Bigl(-\int_\tau^T c\bigl(r, X^{\tau,X(\tau)}(r)\bigr)dr\Bigr) \biggr) \Biggr)\\
=\, & \ee_{t,x}\biggl(\exp\Bigl(-\int_t^\tau c\bigl(r, X^{t,x}(r)\bigr)dr\Bigr) u_1\bigl(\tau, X^{t,x}(\tau)\bigr)\biggr),
\end{align*}
where $u_1 (t, x )$ corresponds to the first term in formula \eqref{eq_143}
\[
u_1(t,x) = \ee_{t,x}\biggl(h\bigl(X^{t,x}(T)\bigr)\exp\Bigl(-\int_t^T c\bigl(r, X^{t,x}(r)\bigr)dr\Bigr) \biggr).
\]

For the second integral in \eqref{eq_148}, we obtain
\begin{align*}
\ee_{t,x}\Biggl(&\ee\biggl(\int_t^T f\bigl(s,X^{t,x}(s)\bigr)\exp\Bigl(-\int_t^s c\bigl(r, X^{t,x}(r)\bigr)dr\Bigr)ds\bigg|\ef^t_\tau\biggr)\Biggr)\\
=\, &\ee_{t,x}\Biggl(\ee\biggl(\int_t^\tau f\bigl(s,X^{t,x}(s)\bigr)\exp\Bigl(-\int_t^s c\bigl(r, X^{t,x}(r)\bigr)dr\Bigr)ds\bigg|\ef^t_\tau\biggr)\Biggr)\\
&+\, \ee_{t,x}\Biggl(\ee\biggl(\int_\tau^T f\bigl(s,X^{t,x}(s)\bigr)\exp\Bigl(-\int_t^s c\bigl(r, X^{t,x}(r)\bigr)dr\Bigr)ds\bigg|\ef^t_\tau\biggr)\Biggr)\\
=\, &\ee_{t,x}\biggl(\int_t^\tau f\bigl(s,X^{t,x}(s)\bigr)\exp\Bigl(-\int_t^s c\bigl(r, X^{t,x}(r)\bigr)dr\Bigr)ds\biggr)\\
&+\, \ee_{t,x}\Biggl(\ee\biggl(\int_\tau^T f\bigl(s,X^{t,x}(s)\bigr)\exp\Bigl(-\int_t^s c\bigl(r, X^{t,x}(r)\bigr)dr\Bigr)ds\bigg|\ef^t_\tau\biggr)\Biggr).
\end{align*}
Making a change of variables in the last integral and using the Markov property of $X(s)$, we get
\begin{align*}
\ee_{t,x}\Biggl(&\ee\biggl(\int_\tau^T f\bigl(s,X^{t,x}(s)\bigr)\exp\Bigl(-\int_t^s c\bigl(r, X^{t,x}(r)\bigr)dr\Bigr)ds\bigg|\ef^t_\tau\biggr)\Biggr)\\
=\, & \ee_{t,x}\Biggl(\exp\Bigl(-\int_t^\tau c\bigl(r, X^{t,x}(r)\bigr)dr\Bigr)\biggr.\\
\biggl.&\times \ee\biggl(\int_\tau^T f\bigl(s,X^{t,x}(s)\bigr)\exp\Bigl(-\int_\tau^s c\bigl(r, X^{t,x}(r)\bigr)dr\Bigr)ds\bigg|\ef^t_\tau\biggr)\Biggr)\\
=\, & \ee_{t,x}\Biggl(\exp\Bigl(-\int_t^\tau c\bigl(r, X^{t,x}(r)\bigr)dr\Bigr)\biggr.\\
\biggl.&\times \ee_{\tau,X(\tau)}\biggl(\int_\tau^T f\bigl(s,X^{\tau,X(\tau)}(s)\bigr)\exp\Bigl(-\int_\tau^s c\bigl(r, X^{\tau,X(\tau)}(r)\bigr)dr\Bigr)ds\biggr)\Biggr)\\
=\, & \ee_{t,x}\biggl(\exp\Bigl(-\int_t^\tau c\bigl(r, X^{t,x}(r)\bigr)dr\Bigr)u_2\bigl(\tau, X^{t,x}(\tau)\bigr)\biggr),
\end{align*}
where $u_2 (t, x )$ corresponds to the second term in formula \eqref{eq_143}
\[
u_2(t,x) = \ee_{t,x}\biggl(\int_t^T f\bigl(s,X^{t,x}(s)\bigr)\exp\Bigl(-\int_t^s c\bigl(r, X^{t,x}(r)\bigr)dr\Bigr)ds \biggr).
\]

Collecting all terms we have
\begin{align*}
w(t,x) = \, &\ee_{t,x}\biggl(u\bigl((\tau,X^{t,x}(\tau)\bigr)\exp\Bigl(-\int_t^\tau c\bigl(r, X^{t,x}(r)\bigr)dr\Bigr)\biggr)\\
& + \ee_{t,x}\biggl( \int_t^\tau f\bigl(s,X^{t,x}(s)\bigr)\exp\Bigl(-\int_t^s c\bigl(r, X^{t,x}(r)\bigr)dr\Bigr)ds\biggr)
\end{align*}
which gives the equality $w(t,x) = v(t,x)$.

Therefore, $w(t,\cdot)\in H^{2+\beta}( U)$ and $u(t,\cdot)\in H^{2+ \beta}_{\mathit{loc}} (\er^d)$ for $t\in (0,T)$ which by Theorem 16.1 of \cite{L1992} implies $u\in C^{1, 2 }\bigl((0,T)\times \er^d\bigr)$ and shows that $u(t,x)$ is a unique classical solution of the Cauchy problem \eqref{eq_142}.
\qed

\medskip

The result of Theorem \ref{th_142} can be easily extended to a theorem showing that a solution of problem \eqref{eq_142} is as smooth as the coefficients of the equation. The theorem below can be obtained from Theorem 16.1 of \cite{L1992} after applying the approximation of a continuous initial condition by a smooth initial condition in the spirit of Lemma 2.6 in \cite{JS2006}.

\bt\label{th_143}
Let $U$ be an arbitrary open bounded set $U$ in $\er^d$ with $\bar U \subset B(R)$. If $h\in BC(U)$, $c(t,x) \ge c_0>0$ for $(t,x)\in (0,T)\times U$, and the coefficients $a_{ij}$, $b_i$, $c$, and $f$ belong to $H^{p+\beta} ( U )$, $\beta\in(0,1)$, with
\begin{align*}
\|a_{ij}(t, \cdot)\|_{H^{p+\beta}(U)} + \|b_i(t, \cdot)&\|_{H^{p+\beta}(U)} +\|c(t, \cdot)\|_{H^{p+\beta}(U)} + \|f(t, \cdot)\|_{H^{p+\beta}(U)}\\
&\le C(T,R), \quad t\in (0,T),
\end{align*}
and $u(t,x)$ is a solution of equation \eqref{eq_142} in $(0,T)\times U$, then $u(t,\cdot)\in H^{p+2+ \beta} (U)$.
\et
\medskip

\begin{remark}
The condition $c_0 > 0$ in the growth assumption of Theorem \ref{th_141}, important for the maximum principle used in the proofs of Theorems \ref{th_142} and \ref{th_143}, can be relaxed to $c_0\in \er$ by a change of variables. If $u(t,x)$ is a solution of the Kolmogorov equation \eqref{eq_142} with $c(t,x) \ge c_0>0$, then $v(t,x) = \e^{-\gamma t}u(t,x)$ solves problem \eqref{eq_142} with $c'(t,x) = c(t,x) -\gamma$, $f'(t,x) = \e^{\gamma t}f(t,x)$. Hence, we can replace $c_0$ with $c'_0 = c_0 -\gamma$ where $c'_0$ can be an arbitrary real number.
\end{remark}

\medskip

\bt\label{th_144} %
Let for an integer $p \ge 2$ the coefficients of equation \eqref{eq_141} fulfill the conditions:
\begin{enumerate}
\item[( S1)] $|\bb(t,\bbx)| +|\bs(t,\bbx)| \le C(T) (1+|\bbx|)$ and $\bb(t,\bbx)$, $\bs(t,\bbx)$ are $(t,\bbx)$-continuous for $(t,x) \in(0,T)\times\er^d$.
\item[(S2)] The derivatives $D_x^\alpha \bb(t,\bbx)$, $D_x^\alpha \bs(t,\bbx)$ exist for $|\alpha|\le p$ and are continuous with respect to $\bbx$ for all $t\in (0,T)$, $\bbx\in \er^d$.
\item[(S3)] $\big|D_x^\alpha \bb(t,\bbx)\big| + \big|D_x^\alpha \bs(t,\bbx)\big| \le C$ for $0<|\alpha|\le p$ and $(t,x) \in(0,T)\times\er^d$.
\end{enumerate}

Assume that for the same $p\ge 2$ the coefficients of the Cauchy problem \eqref{eq_142} are such that:
\begin{enumerate}
\item[( K1)] There exists $\delta(T) >0$ such that for $(t,x) \in [0,T]\times\er^d$ and $\xi\in \er^d\setminus \{0\}$
\[
\sum_{i,j =1}^d a_{ij}(t,x) \xi_i\xi_j \ge \delta(T) |\xi|^2.
\]
\item[( K2)] Functions $c\sr[0,\infty) \times \er^d\to [0,\infty)$ and $f\sr[0,\infty) \times \er^d\to \er$ are $(t,x)$-continuous and $c(t,\bbx)\ge c_0>0$. Function $h\sr \er^d \to \er$ belongs to $BC(\er^d)$.
\item[( K3)] The derivatives $D_x^\alpha c(t,\bbx)$, $D_x^\alpha f(t,\bbx)$, $D_x^\alpha h( \bbx)$ exist for $|\alpha|\le p$ and are continuous with respect to $\bbx$ for all $t\in (0,T)$.
\item[(K4)] For all $t\in (0,T)$ and $|\alpha|\le p$ there exist $C>0$, $q\ge 1$ such that
\[
\big|D_x^\alpha c(t,\bbx)\big| + \big|D_x^\alpha f(t,\bbx)\big| + \big|D_x^\alpha h( \bbx)\big| \le C(1 + |\bbx|^{2q}), \quad x\in\er^d.
\]
\end{enumerate}
Then $u(t,x)$ given by equation \eqref{eq_143} is a well defined function in $[0,T]\times\er^d$ which is continuously differentiable with respect to $t$ in $(0,T)\times\er^d$ and has $(t,\bbx)$-continuous derivatives $D_x^\alpha u(t,\bbx)$, $|\alpha| \le p$, with the growth estimate
\begin{equation}\label{eq_148a}
\sup_{0< t< T} \big|D_x^\alpha u(t,\bbx)\big| \le C(1 +|\bbx|^{2q(1+|\alpha|)}), \quad |\alpha|\le p,\ \bbx\in \er^d,
\end{equation}
where $C$ is a positive constant, and $q$ is the same as in the growth estimates of the derivatives of $c$, $f$, and $h$.
\et
\proof
The assertion that $u(t,x)$ is continuously differentiable with respect to $t$ and possesses continuous derivatives $D_x^\alpha u(t,\bbx)$, $|\alpha| \le p$, follows straightforwardly from Theorem \ref{th_143}. What we want to prove is the possibility of differentiation under the integral sign. This differentiation is also crucial to show the estimates of the derivatives.

To simplify the proof, we take $d =1$. First, we will show that under the assumption of the theorem $u(t,x)$ given by formula \eqref{eq_143} can be continuously differentiated with respect to $\bbx$ under the integral sign and $D^p_xu(t,x) \le C(1 +|\bbx|^{2q(1+p)})$.

Consider the difference $u(t, x +\Delta x) - u(t,x)$. We can write this difference as the sum of terms:
\begin{align*}
u(&t,x +\Delta x) - u(t,x)\\
= &\, \ee\biggl(\Bigl( h\bigl(X^{t,x+\Delta x}(T)\bigr) - h\bigl(X^{t,x}(T)\bigr)\Bigr)\exp\Bigl(-\int_t^T c\bigl(r, X^{t,x}(r)\bigr)dr\Bigr)\biggr)\\
&+ \ee\Biggl(h\bigl(X^{t,x+\Delta x}(T)\bigr)\Biggr.\\
\Biggl. &\times\biggl(\exp\Bigl(-\int_t^T c\bigl(r, X^{t,x+\Delta x}(r)\bigr)dr\Bigr) - \exp\Bigl(-\int_t^T c\bigl(r, X^{t,x}(r)\bigr)dr\Bigr)\biggr)\Biggr)\\
&+ \ee \biggl(\int_t^T \Bigl(f\bigl(s,X^{t,x+\Delta x}(s)\bigr) -f\bigl(s,X^{t,x}(s)\bigr)\Bigr)\biggr.\\
\biggl. &\phantom{f\bigl(s,X^{t,x}(s)\bigr) f\bigl(s, (s)\bigr) }\times\exp\Bigl(-\int_t^s c\bigl(r, X^{t,x}(r)\bigr)dr\Bigr)ds\biggr)\\
&+ \ee \Biggl(\int_t^T f\bigl(s,X^{t,x+\Delta x}(s)\bigr) \Biggr.\\
\Biggl. &\times\biggl(\exp\Bigl(-\int_t^s c\bigl(r, X^{t,x+\Delta x}(r)\bigr)dr\Bigr) - \exp\Bigl(-\int_t^s c\bigl(r, X^{t,x}(r)\bigr)dr\Bigr)\biggr)ds\Biggr).
\end{align*}

Our goal is to show that there exists a limit $\bigl(u(t,x +\Delta x) - u(t,x)\bigr)/\Delta x$ for $\Delta x \to 0$.
We will consider the convergence for each term of the difference $\bigl(u(t,x +\Delta x) - u(t,x)\bigr)$ separately. For the first term, we have by the mean value theorem\\
\begin{align*}
\ee\biggl(\Bigl(& h\bigl(X^{t,x+\Delta x}(T)\bigr) - h\bigl(X^{t,x}(T)\bigr)\Bigr)\exp\Bigl(-\int_t^T c\bigl(r, X^{t,x}(r)\bigr)dr\Bigr)\biggr)\\
=\, &\ee\biggl(\int_0^1 \frac{d}{d\mu} h\Bigl(X^{t,x}(T) + \mu\bigl(X^{t,x+\Delta x}(T) -X^{t,x}(T)\bigr)\Bigr)d\mu \biggr.\\
\biggl. &\times\exp\Bigl(-\int_t^T c\bigl(r, X^{t,x}(r)\bigr)dr\Bigr)\biggr)\\
=\, &\ee\biggl(\int_0^1 D_x h\Bigl(X^{t,x}(T) + \mu\bigl(X^{t,x+\Delta x}(T) -X^{t,x}(T)\bigr)\Bigr)d\mu \biggr.\\
\biggl. &\times \bigl(X^{t,x+\Delta x}(T) -X^{t,x}(T)\bigr) \exp\Bigl(-\int_t^T c\bigl(r, X^{t,x}(r)\bigr)dr\Bigr)\biggr).
\end{align*}
Dividing the above expression by $\Delta x$, we get for $\Delta x\to 0$ the convergence
\[
\frac{X^{t,x+\Delta x}(T) -(X^{t,x}(T)}{\Delta x}\to \frac{\partial}{\partial x} X^{t,x}(T)
\]
in $L^2$ due to Theorem 8.5.1 in \cite{GS1969}. We have also for $\Delta x\to 0$
\[
\int_0^1 D_xh\Bigl(X^{t,x}(T) + \mu\bigl(X^{t,x+\Delta x}(T) -X^{t,x}(T)\bigr)\Bigr)d\mu \to D_xh\bigl(X^{t,x}(T)\bigr)
\]
in $L^2$. This last convergence follows from the $x$-continuity of $X^{t,x}(s)$, Lemma III.6.13 in \cite {K1995}, and the estimate
\[
\ee\biggl(\int_0^1 \!\! D_xh\Bigl(X^{t,x}(T) + \mu\bigl(X^{t,x+\Delta x}(T) -X^{t,x}(T)\bigr)\Bigr)d\mu - D_xh \bigl((X^{t,x}(T)\bigr)\biggr)^4\le C,
\]
which holds for any $|x|\le C$ due to Theorem II.2.1 in \cite{K1984} and the growth condition for $D_xh(x)$.

Hence, due to the continuity of $c(t,x)$, the expression
\[
\frac{1}{\Delta x}\ee\biggl(\Bigl( h\bigl(X^{t,x+\Delta x}(T)\bigr) - h\bigl(X^{t,x}(T)\bigr)\Bigr)\exp\Bigl(-\int_t^T c\bigl(r, X^{t,x}(r)\bigr)dr\Bigr)\biggr)
\]
converges to
\[
\ee\biggl( D_x h\bigl(X^{t,x}(T)\bigr)\frac{\partial}{\partial x} X^{t,x}(T) \exp\Bigl(-\int_t^T c\bigl(r, X^{t,x}(r)\bigr)dr\Bigr)\biggr).
\]

Similar considerations show that the term
\begin{align*}
\frac{1}{\Delta x} \ee &\biggl(\int_t^T \Bigl(f\bigl(s,X^{t,x+\Delta x}(s)\bigr) -f\bigl(s,X^{t,x}(s)\bigr)\Bigr)\biggr.\\
\biggl. &\phantom{f\bigl(s,X^{t,x}(s)\bigr) f\bigl(s,X^{t,x}(s)\bigr) }\times\exp\Bigl(-\int_t^s c\bigl(r, X^{t,x}(r)\bigr)dr\Bigr)ds\biggr)
\end{align*}
converges in $L^2$ to
\[
\ee \biggl(\int_t^T D_x f\bigl(s,X^{t,x }(s)\bigr) \frac{\partial}{\partial x} X^{t,x}(s)\exp\Bigl(-\int_t^s c\bigl(r, X^{t,x}(r)\bigr)dr\Bigr)ds\biggr).
\]

To obtain convergence for the two other terms let us observe the equality
\begin{align*}
\exp\Bigl(&-\int_t^s c\bigl(r, X^{t,x+\Delta x}(r)\bigr)dr\Bigr) - \exp\Bigl(-\int_t^s c\bigl(r, X^{t,x}(r)\bigr)dr\Bigr)\\
=\, &\int_0^1 \frac{d}{d\mu} \exp\biggl(-\int_t^s c\Bigl(r, X^{t,x}(r) +\mu\bigl( X^{t,x+\Delta x}(r) - X^{t,x}(r)\bigr)\Bigr)dr\biggr)d\mu\\
=\, &- \int_0^1 \Biggl( \exp\biggl(-\int_t^s c\Bigl(r, X^{t,x}(r) +\mu\bigl( X^{t,x+\Delta x}(r) - X^{t,x}(r)\bigr)\Bigr)dr\biggr)\Biggr.\\
&\phantom{ X^{t,x+\Delta x}(r)} \times \int_t^s \biggl(D_x c\Bigl(r, X^{t,x}(r) +\mu\bigl( X^{t,x+\Delta x}(r) - X^{t,x}(r)\bigr)\Bigr)\biggr.\\
&\Biggl.\biggl. \phantom{ X^{t,x+\Delta x}(r) - X^{t,x}(r) }\times \bigl( X^{t,x+\Delta x}(r) - X^{t,x}(r)\bigr)\biggr) dr\Biggr) d\mu.
\end{align*}

By Lemma III.6.13 in \cite {K1995}, Theorem II.2.1 in \cite{K1984}, the continuity of $h(x)$, $f(t,x)$ and $c(t,x)$, and the growth conditions for $h(x)$, $f(t,x)$, $c(t,x)$, and $D_x c(t,x)$, we prove that
\begin{align*}
\frac{1}{\Delta x}\ee&\Biggl(h\bigl(X^{t,x+\Delta x}(T)\bigr)\biggr.\\
\biggl. &\times\biggl(\exp\Bigl(-\int_t^T c\bigl(r, X^{t,x+\Delta x}(r)\bigr)dr\Bigr) - \exp\Bigl(-\int_t^T c\bigl(r, X^{t,x}(r)\bigr)dr\Bigr)\biggr)\Biggr)
\end{align*}
converges to
\begin{align*}
- \ee\biggl(h\bigl(X^{t,x}(T)\bigr)
\exp\Bigl(-\int_t^T \!\!\! c\bigl(r, X^{t,x}(r)\bigr)dr\Bigr)
\!\!\int_t^T \!\!\! D_xc\bigl(r, X^{t,x}(r)\bigr)\frac{\partial}{\partial x} X^{t,x}(r) dr
\biggr),
\end{align*}
and
\begin{align*}
\frac{1}{\Delta x}&\ee \Biggl(\int_t^T f\bigl(s,X^{t,x+\Delta x}(s)\bigr) \biggr.\\
\biggl. &\times\biggl(\exp\Bigl(-\int_t^s c\bigl(r, X^{t,x+\Delta x}(r)\bigr)dr\Bigr) - \exp\Bigl(-\int_t^s c\bigl(r, X^{t,x}(r)\bigr)dr\Bigr)\biggr)ds\Biggr)
\end{align*}
converges to
\begin{align*}
- \ee\Biggl(\int_t^T \biggl(f\bigl(X^{t,x}(s)\bigr)
&\exp\Bigl(-\int_t^s c\bigl(r, X^{t,x}(r)\bigr)dr\Bigr)\biggr. \\
&\biggl.\times\int_t^s D_x c\bigl(r, X^{t,x}(r)\bigr)\frac{\partial}{\partial x} X^{t,x}(r) dr\biggr)ds \Biggr).
\end{align*}
Together the above expressions prove that we can differentiate \eqref{eq_143} under the integral sign and obtain the following formula for the $x$-derivative of $u(t,x)$
\begin{align*}
D_x &u(t,x) =\, \ee\biggl(D_x h\bigl(X^{t,x}(T)\bigr)\frac{\partial}{\partial x} X^{t,x}(T) \exp\Bigl(-\int_t^T c\bigl(r, X^{t,x}(r)\bigr)dr\Bigr)\biggr)\\
& - \ee\biggl(h\bigl(X^{t,x}(T)\bigr)
\exp\Bigl(-\int_t^T c\bigl(r, X^{t,x}(r)\bigr)dr\Bigr)\biggr.\\
&\biggl.\phantom{c\bigl(r, X^{t,x}(r)\bigr) } \times \int_t^T D_x c\bigl(r, X^{t,x}(r)\bigr)\frac{\partial}{\partial x} X^{t,x}(r) dr
\biggr)\\
& + \ee \biggl(\int_t^T D_x f\bigl(s,X^{t,x }(s)\bigr) \frac{\partial}{\partial x} X^{t,x}(s)\exp\Bigl(-\int_t^s c\bigl(r, X^{t,x}(r)\bigr)dr\Bigr)ds\biggr)\\
& - \ee\Biggl(\int_t^T \biggl(f\bigl(X^{t,x}(s)\bigr)
\exp\Bigl(-\int_t^s c\bigl(r, X^{t,x}(r)\bigr)dr\Bigr)\biggr. \\
&\biggl. \phantom{c\bigl(r, X^{t,x}(r)\bigr) }\times\int_t^s D_x c\bigl(r, X^{t,x}(r)\bigr)\frac{\partial}{\partial x} X^{t,x}(r) dr\biggr)ds \Biggr).
\end{align*}

The estimate of $D_x u(t,x)$ can be obtain similarly like in the proof of Theorem \ref{th_141} taking into account the above expression for $D_x u(t,x)$ and the boundedness of $|D_x X^{t,x}(s)|$
\begin{align*}
|D_x u(t,x)| \le &\, C\, \ee\Big|D_x h\bigl(X^{t,x}(T)\bigr) \Big| \\
&+ C\, \ee\biggl(\big|h\bigl(X^{t,x}(T)\bigr)\big|\int_t^T \big| D_x c\bigl(r, X^{t,x}(r)\bigr)\big| dr \biggr)\\
& + C\, \ee \biggl(\int_t^T \big|D_x f\bigl(s,X^{t,x }(s)\bigr)\big| ds\biggr)\\
& +C\, \ee\Biggl(\int_t^T \biggl(\big|f\bigl(X^{t,x}(s)\bigr)\big|
\int_t^s \big|D_x c\bigl(r, X^{t,x}(r)\bigr)\big|dr \biggr)ds \Biggr).
\end{align*}
Straightforwardly we obtain
\begin{align*}
\ee\Big|D_x h\bigl(X^{t,x}(T)\bigr) \Big| + \ee \biggl(\int_t^T \big|D_x f\bigl(s,X^{t,x }(s)\bigr)\big| ds\biggr) \le C\, \bigl(1 + |x|^{2q}\bigr).
\end{align*}
To estimate the other two terms, we apply H\"older's inequality and get
\begin{align*}
&\ee\biggl(\big|h\bigl(X^{t,x}(T)\bigr)\big|\int_t^T \big| D_x c\bigl(r, X^{t,x}(r)\bigr)\big| dr \biggr) \\
& \phantom{h\bigl(X^{t,x}(T)\bigr)} \le C\, \Bigl(\ee \big|h\bigl(X^{t,x}(T)\bigr)\big|^2\Bigr)^{1/2} \biggl(\int_t^T \ee\big| D_x c\bigl(r, X^{t,x}(r)\bigr)\big|^2 dr \biggr)^{1/2}\\
& \phantom{h\bigl(X^{t,x}(T)\bigr)} \le C\, \bigl(1+|x|^{4q}\bigr),
\end{align*}
\begin{align*}
&\ee\Biggl(\int_t^T \biggl(\big|f\bigl(X^{t,x}(s)\bigr)\big|
\int_t^s \big|D_x c\bigl(r, X^{t,x}(r)\bigr)\big|dr \biggr)ds \Biggr) \\
& \phantom{h\bigl(X^{t,x}(T)\bigr)} \le C\, \int_t^T \Bigl(\ee \big|f\bigl(X^{t,x}(s)\bigr)\big|^2\Bigr)^{1/2} \biggl(\int_t^s \ee\big| D_x c\bigl(r, X^{t,x}(r)\bigr)\big|^2 dr \biggr)^{1/2} ds\\
& \phantom{h\bigl(X^{t,x}(T)\bigr)} \le C\, \bigl(1+|x|^{4q}\bigr).
\end{align*}

For higher $x$-derivatives of $u(t,x)$ we can proceed analogously.
\qed

The following corollary reformulates the results of Theorem \ref{th_144} in terms of the partial Schauder estimates.

\begin{cor}\label{cor_141}
Let the assumptions of Theorem \ref{th_144} hold for $p\ge 2$. We assume additionally that there exists a constant $C>0$ such that for all $t\in [0,T]$
\[
|D^\alpha_x c(t,x)| \le C,\ 0<|\alpha|\le p, \quad x\in \er^d.
\]
Then for the function $u(t,x)$ defined by equation \eqref{eq_143} the following partial Schau\-der estimates in the space of differentiable functions with the polynomial weight $P(x) = 1+|x|^{2q}$ are valid
\[
\sup_{0<t<T}\|u(t,\cdot)\|_{BC^p_P(\er^d)} \le C \Bigl(\|h\|_{BC^p_P(\er^d)} + \sup_{0<t<T}\|f(t,\cdot)\|_{BC^p_P(\er^d)}\Bigr).
\]
\end{cor}
\proof
The estimate in $BC_P(\er^d)$ can be obtained by repeating the computations in the proof of Theorem \ref{th_141}. Let us recall that $u(t,x)$ is given by the Feynman-Kac formula
\begin{equation*}
\begin{split}
u(t,x) = \, &\ee\biggl(h\bigl(X^{t,x}(T)\bigr)\exp\Bigl(-\int_t^T c\bigl(r, X^{t,x}(r)\bigr)dr\Bigr)\biggr.\\
& \biggl.+ \int_t^T f\bigl(s,X^{t,x}(s)\bigr)\exp\Bigl(-\int_t^s c\bigl(r, X^{t,x}(r)\bigr)dr\Bigr)ds\biggr).
\end{split}
\end{equation*}
For the first term of this formula, we have the estimate
\begin{align*}
\sup_{0<t<T}&\sup_{x\in\er^d} \frac1{P(x)}\bigg| \ee\biggl( h\bigl(X^{t,x}(T)\bigr) \exp\Bigl(-\int_t^T c\bigl(r, X^{t,x}(r)\bigr)dr\Bigr)\biggr)\bigg|\\
&\le \sup_{0<t<T}\sup_{x\in\er^d} \frac1{P(x)} \ee \Big|h\bigl(X^{t,x}(T)\bigr)\Big|\\
&= \sup_{0<t<T}\sup_{x\in\er^d} \frac1{P(x)} \ee \Bigg|\frac{h\bigl(X^{t,x}(T)\bigr)}{1 +\bigl|X^{t,x}(T)\bigr|^{2q}} \Bigl(1 + \bigl|X^{t,x}(T)\bigr|^{2q} \Bigr) \Bigg| \\
&\le C \|h\|_{BC(\er^d)} \sup_{0<t<T}\sup_{x\in\er^d} \frac1{P(x)} \Bigl(1 + \ee\bigl|X^{t,x}(T)\bigr|^{2q} \Bigr)
\le C\|h\|_{BC(\er^d)}.
\end{align*}

For the second term in the Feynman-Kac formula we obtain
\begin{align*}
\sup_{0<t<T}&\sup_{x\in\er^d} \frac1{P(x)}\bigg|\ee\biggl(\int_{t}^T f\bigl(s, X^{t,x}(s)\bigr) \exp\Bigl(-\int_{t}^s c\bigl(r, X^{t,x}(r)\bigr)dr\Bigr)ds\biggr)\bigg|\\
&\le \, C \sup_{0<t<T}\sup_{x\in\er^d} \frac1{P(x)}\ee \biggl(\int_{t}^T \Big|f\bigl(s, X^{t,x}(s)\bigr)\Big| ds\biggr)\\
&= C \sup_{0<t<T}\sup_{x\in\er^d} \frac1{P(x)}\ee \biggl(\int_{t}^T \frac{\Big|f\bigl(s, X^{t,x}(s)\bigr)\Big|}{1 +\bigl|X^{t,x}(s)\bigr|^{2q}} \Bigl(1 + \bigl|X^{t,x}(s)\bigr|^{2q} \Bigr) ds\biggr)\\
&\le \, C \ \sup_{0<t<T} \|f(t,\cdot)\|_{BC(\er^d)}\sup_{x\in\er^d} \frac1{P(x)} \ee \biggl(\int_{t }^T \Bigl(1 + \big|X^{t ,x }(s)\big|^{2q}\Bigr)ds\biggr) \\
&\le \, C \sup_{0<t<T}\|f(t,\cdot)\|_{BC(\er^d)}.
\end{align*}
Together these estimates give the desired Schauder estimate in $BC_P(\er^d)$.

The estimates in $BC_P^p(\er^d)$ for $p>0$ can be obtained similarly. By assumption, the derivatives $D^\alpha_x c(t,x)$ and $D^\alpha_x X^{t,x}(s)$ for $0<|\alpha|\le p$ are bounded, then the derivative of $u(t,x)$ with respect to $x_j$ has the estimate
\begin{align*}
|D^j_x u(t,x)| \le& \, C\, \ee\Big|D^j_x h\bigl(X^{t,x}(T)\bigr) \Big|
+ C\, \ee \big|h\bigl(X^{t,x}(T)\bigr)\big| \\
& + C\, \ee \biggl(\int_t^T \!\!\big|D^j_x f\bigl(s,X^{t,x }(s)\bigr)\big| ds\biggr)
+C\, \ee\biggl(\int_t^T \!\!\big|f\bigl(X^{t,x}(s)\bigr)\big|
ds \biggr).
\end{align*}
Performing similar computations as for $u(t,x)$ and applying Lemma 8.4.1 from \cite{GS1969}, we obtain
\[
\|D^j_x u(t,\cdot)\|_{BC_P(\er^d)} \le C \| h \|_{BC^1_P(\er^d)} + C \sup_{0<t<T}\|f(t,\cdot)\|_{BC^1_P(\er^d)}.
\]

Computations for higher derivatives are analogous.
\qed

\section{Schauder's estimates}

This section aims to give an analytic proof of the partial Schauder estimates similar to the estimates of Corollary \ref{cor_141} omitting the stochastic representation of the solution given by the Feynman-Kac formula.
The partial Schauder estimates for the Kolmogorov equation in weighted H\"ol\-der's spaces have been obtained by Lorenzi in \cite{L2000} under the assumption of bound\-ed coefficients. Here, we prove the extension of his result to unbounded coefficients. The first step in this approach is to obtain known Schauder's estimates in spaces of bounded H\"older's functions. We follow in this presentation the approach from \cite{L2011} (see also \cite{BL2005} where similar methods were used). There is a non obvious choice in proving the partial Schauder estimates concerning the assumed time regularity of the coefficients. In already mentioned papers \cite{L2000} and \cite{L2011}, the author assumes only time measurability of the coefficients (at least for a part of the results). We assume that the coefficients are continuous for $t\in [0,T]$ ($T$ finite). This assumption is motivated by our consideration of classical solutions of the Kolmogorov equation, i.e., solutions that are of class $C^{1,2}$. Thus, we get strong stochastic solutions by standard assumptions of $(t,x)$-continuity of the coefficients in the stochastic equation and classical solutions of the Kolmogorov equation assuming $t$-continuity and $x$-H\"older continuity of its coefficients.

We begin by proving the existence of solutions for the homogeneous version of the Kolmogorov equation
\begin{equation}\label{eq_161}
\begin{split}
& D_t u(t,x) + \A^tu(t,x) - c(t,x)u(t,x) = 0, \ (t,x) \in [0,T)\times\er^d,\\
&u(T,x) = h(x), \quad x\in \er^d,
\end{split}
\end{equation}
where operator $\A^t$ and functions $c(t,x)$ and $h(x)$ are the same as in equation \eqref{eq_142}.

From Theorem \ref{th_142} we can easily obtain the following result.

\bt\label{th_161} %
Let the coefficients of equation \eqref{eq_141} fulfill the conditions of Assumption \ref{as_121} and operator $\A^t$ be uniformly elliptic, i.e., there exists $\delta(T) >0$ such that for $(t,x) \in [0,T]\times\er^d$ and $\xi\in \er^d\setminus \{0\}$
\[
\sum_{i,j =1}^d a_{ij}(t,x) \xi_i\xi_j \ge \delta(T) |\xi|^2.
\]
On the coefficients of equation \eqref{eq_161}, we assume:
\begin{enumerate}
\item Function $c\sr[0,T] \times \er^d\to \er$\, is $(t,x)$-continuous with $c(t,x) \ge c_0 >0$. In addition, $c$ as a function of $x$ is locally H\"older continuous, i.e., for an open bounded set $U$ in $\er^d$ such that $\bar U \subset B(R)$, $c(t,\cdot) \in H^{\beta} (U )$, $0< \beta<1$, and $\|c(t, \cdot)\|_{H^{\beta}(U)} \le C(T,R)$.
\item Function $h\sr \er^d \to \er$\, belongs to $BC(\er^d)$.
\end{enumerate}
Under the above assumptions the function $u(t,x)$ given by the equation
\begin{equation}\label{eq_162}
u(t,x) = \, \ee\biggl(h\bigl(X^{t,x}(T)\bigr)\exp\Bigl(-\int_t^T c\bigl(r, X^{t,x}(r)\bigr)dr\Bigr)\biggr)
\end{equation}
is such that $u(t,\cdot) \in H^{2+\beta} (U)$ and is a bounded function in $[0,T]\times\er^d$ with
\[
\sup_{[0,T]\times\er^d} |u(t,x| \le C(T) \|h\|_{BC(\er^d)}.
\]
Furthermore, $u(t,x)$ belongs to $C^{1,2}\bigl((0,T)\times\er^d\bigr)$ and is a unique classical solution of the Cauchy problem \eqref{eq_161}.
\et

To obtain Schauder's estimates for a solution of equation \eqref{eq_161} we need a maxi\-mum principle.
The following maximum principle is a simple corollary from Proposition 2.1 in \cite{L1998}.

\bt\label{th_162}
Let $u \in BC([0, T ] \times \er^d)\cap C^{1, 2} ([0, T)\times \er^d)$ be a solution of the Cauchy problem
\begin{equation*}
\begin{split}
& D_t u(t,x) + \A^tu(t,x) - c(t,x)u(t,x) + f(t,x)=0, \\
&u(T,x) = h(x),
\end{split}
\end{equation*}
with $c(t,x) \ge c_0$.

If $f(t,x) \le 0$ for $(t,x) \in [0, T ] \times \er^d$ and $\sup_{[0, T ] \times \er^d} u(t,x) > 0$, then
\[
\sup_{x\in \er^d} u(t,x) \le \e^{-c_0(T-t)} \sup_{x\in \er^d} h(x), \quad t\in [0,T].
\]

If $f(t,x) \ge 0$ for $(t,x) \in [0, T ] \times \er^d$ and $\inf_{[0, T ] \times \er^d} u(t,x) < 0$, then
\[
\inf_{x\in \er^d} u(t,x) \ge \e^{-c_0(T-t)} \inf_{x\in \er^d} h(x), \quad t\in [0,T].
\]

In particular, if $f(t,x)=0$ then
\[
\|u(t, \cdot)\|_{BC(\er^d)} \le \e^{-c_0(T-t)} \|h\|_{BC(\er^d)}, \quad t\in [0,T].
\]
\et
\proof
Observe that in the proof of Proposition 2.1 in \cite{L1998} no regularity of coefficients is necessary. The proof requires only the uniform ellipticity of $\A^t$, the lower bound on $c(t,x)$, and the existence of the Lyapunov function $\phi\in C^2(\er^d)$ such that for a constant $\lambda >0$
\[
\lim_{|x|\to +\infty} \phi(x) = +\infty,\quad \A^t\phi-c(t,x)\phi \le \lambda\phi, \quad (t,x)\in [0,T]\times\er^d.
\]
Taking into account the growth conditions for $\bb(t,x)$ and $\bs(t,x)$, we can take $\phi(x) = (1+|x|^2)$. Then we get
\[
\A^t \phi - c(t,x)\phi \le \big|\A^t (1+|x|^2)\big| - c_0(1+|x|^2)\le
C(1+|x|)^2 + |c_0|(1+|x|^2) \le \lambda\phi,
\]
where $\lambda > |c_0| + C$.
\qed

\begin{defi}
Let $u(t,x)$ be a solution of problem \eqref{eq_161} with initial data $u(s,x) = h(x)$, $s>t$. $G(\cdot, s)$ is defined as an operator on $BC(\er^d)$ which associates with any $h(x)\in BC(\er^d)$ a unique classical solution of problem \eqref{eq_161}.
\end{defi}

The following corollary is a simple conclusion from Theorems \ref{th_161} and \ref{th_162}.

\begin{cor}
For any $s > t$ and any $h \in BC(\er^d)$, define $G ( t, s ) h = u ( t, \cdot)$, where $u ( t, \cdot)$ is a unique solution of problem \eqref{eq_161} such that $u(s,\cdot) = h(\cdot)$. Then, $G ( t, s )$ is a bounded linear operator in $BC(\er^d)$ and the following estimate holds
\[
\|G(t,s)\|_{\mathcal{L}\bigl(BC(\er^d), BC(\er^d)\bigr)} \le \e^{-c_0(s-t)}.
\]
\end{cor}

To obtain the Schauder estimates for solutions of equation \eqref{eq_161}, we need additional assumptions that include conditions on the coefficients of the stochastic equation \eqref{eq_141}.

\begin{assu}\label{as_161}
\n
\begin{enumerate}
\item Let the coefficients $\bb_i(t,\bbx)$ and $\bs_{ij}(t,\bbx)$ of equation \eqref{eq_141} be $(t,x)$-conti\-nu\-ous and with linear growth
\[
|\bb(t,\bbx)| +|\bs(t,\bbx)| \le C (1+|\bbx|), \ t\in [0,T],
\]
and three times continuously differentiable with respect to $\bbx$ with derivatives $D^\alpha_x \bb_i(t,\bbx), D^\alpha_x \bs_{ij}(t,\bbx) \in BC(\er^d)$ for $|\alpha|= 1,2, 3$, $t\in [0,T]$, and $D^\alpha_x \bb_i(t,\bbx), D^\alpha_x \bs_{ij}(t,\bbx) \in H^\beta_{\mathit{loc}}(\er^d)$ for $|\alpha| = 3$.
\item Let the growth of $\bs_{ij}(s,y)$ be limited by the inequalities
\[
C_1 (1+|\bbx|)^{1/2} \le \bs_{ij}(t,\bbx) \le C_2 (1+|\bbx|), \quad t\in[0,T],
\]
where $C_1$ and $C_2$ are positive constants.
\item $c(t,x) \ge c_0$ for $(t,x)\in [0,T]\times\er^d$ with $c_0 > 0$. Function $c(t,x)$ is three times continuously differentiable with respect to $x$ and
\[
| D^\alpha_x c(t,x)| \le C, \quad 1\le |\alpha|\le 3.
\]
\end{enumerate}

\end{assu}

As we mentioned earlier, we prove Schauder's estimates in H\"older's spaces of bounded functions following the approach from \cite{L2011}. An important part of this proof is Theorem 2.4 of \cite{L2011}. The main difference with the original assumptions of this theorem is in the regularity of the coefficients. In Theorem 2.4 of \cite{L2011} the author assumes the H\"older regularity of the coefficients with respect to variables $(t,x)$. In our approach, H\"older's regularity is required only with respect to the spatial variables due to the partial Schauder estimates. 
With respect to variables $(t,x)$, we need only uniform continuity on $[0,T]\times\er^d$. The following lemma shows that Assumption \ref{as_161} implies conditions required
in \cite{L2011} for the proof of Theorem 2.4, and after some natural modifications, we can use that theorem.

\begin{lemma}\label{lem_161}
Under Assumption \ref{as_161} the following reformulation of Hypotheses 2.1 of \cite{L2011} is valid for a universal positive constant $C$:
\begin{enumerate}
\item The coefficients $a_{ij}$, $b_i$, $i, j = 1, \dots ,d$, and $c$ belong to $H^{3+\beta}_{\mathit{loc}}(\er^d)$ for some $\beta \in(0, 1)$.
\item Operator $\A^t$ is uniformly elliptic
\[
\sum_{i,j=1}^d a_{ij}(t,x) \xi_i\xi_j \ge \nu(t,x) |\xi|^2,\ (t,x)\in [0,T]\times\er^d,\ \xi \in \er^d \setminus\{0\},
\]
for some function $\nu(t,x) \ge \nu_0>0$, $(t,x)\in [0,T]\times\er^d$.
\item The following estimates hold
\begin{align*}
&\Big|\sum_{j=1}^d a_{ij}(t,x)x_j\Big| \le C(1+|x|^2)\nu(t,x), \\
&\sum_{j=1}^d a_{jj}(t,x) \le C(1+|x|^2)\nu(t,x),\\
&\sum_{j=1}^d b_{j}(t,x)x_j \le C(1+|x|^2)\nu(t,x).
\end{align*}
\item $|D_x^\alpha a_{ij}(t,x)| \le C \nu(t,x)$ for $|\alpha| = 1, 3$ and
\[
\sum_{i,j,k,l=1}^d D_x^{kl}a_{ij}(t,x)m_{ij}m_{kl} \le C \nu(t,x) \sum_{k,l=1}^d |m_{kl}|^2,
\]
for any $d\times d$ symmetric matrix $M = \{m_{kl}\}$ and $(t,x)\in [0,T]\times\er^d$.
\item For the entries of $b$ we have
\begin{align*}
&\sum_{i,j=1}^d D_x^j b_i(t,x) \xi_i\xi_j \le C \nu(t,x) |\xi|^2,\\
& |D_x^\alpha b_i(t,x)| \le C \nu(t,x), \ |\alpha| = 2,3,\ i = 1,\dots, d.
\end{align*}
\item There exist a positive function $\phi\in C^2(\er^d)$ and a constant $\lambda >0$ such that
\[
\lim_{|x|\to +\infty} \phi(x) = +\infty,\quad \sup_{x\in\er^d} \bigl(\A^t\phi(x)-\lambda\phi(x)\bigr) < +\infty.
\]
\end{enumerate}

\end{lemma}
\proof
Since $a_{ij} = \frac12\sum_{k=1}^m \bs_{ik}\bs_{jk}$ then point 2 of Assumption \ref{as_161} gives $\nu(t,x) = C(1+|x|)$ with $C>0$. Therefore points 1 and 2 of the lemma follow straightforwardly from Assumption \ref{as_161}. Point 3 of the lemma holds true due to the linear growth of $\bb_i$ and $\bs_{ij}$. To obtain point 4, let us observe that $|D_x^\alpha a_{ij}(t,x)| \le C (1+|\bbx|)$, for $|\alpha| \le 3$, due to the linear growth of $\bs_{ij}$ and the boundedness of derivatives $D_x^\alpha \bs_{ij}$. The boundedness of $D_x^\alpha \bb_{i}$ for $|\alpha| \le 3$ implies point 5 of the lemma. Point 6 can be proved similarly like in Theorem \ref{th_162} taking $\phi(x) = (1 +|x|^2)$.
\qed
\medskip

The following theorem is a reformulation of Theorem 2.4 from \cite{L2011} (see also Theorem 3.3 in \cite{BL2005}). We present an outline of the proof pointing to the differences with the original proof implied by the modifications in estimates due to Assumption \ref{as_161}.
As we pointed out before, the main
difference with paper \cite{L2011} is in the regularity of the coefficients. We will show that the H\"older regularity with respect to variables $(t,x)$ can be replaced by a weaker regularity only in the spatial variables and $t$-continuity.

\bt\label{th_163}
Let Assumption \ref{as_161} be satisfied. Then, for any $p_1+\beta_1, p_2+\beta_2\in [0, 3]$, with $p_1+\beta_1 \le p_2+\beta_2$, $\beta_1,\beta_2 \in [0,1)$, there exists a positive constant $C$ such that
\[
\|G(t,T)h\|_{H^{p_2+\beta_2}(\er^d)} \le C(T-t)^{-\frac{p_2+\beta_2 -p_1-\beta_1}{2}}\|h\|_{H^{p_1+\beta_1}(\er^d)},
\]
where $H^{p}(\er^d)$ for an integer $p$ denotes the space $BC^p(\er^d)$.
\et
\proof
We present only the part of the proof that corresponds to $\beta_1=\beta_2 =0$ and $p_1 = 0$, $p_2 =3$. The idea of the proof in \cite{L2011} is to apply Bernstein's method to equation \eqref{eq_142} restricted to an open, bounded ball $B(R)$ with radius $R$. Due to Theorem \ref{th_142} there is a function $u_R(t,x)$ that is a classical solution for the problem
\begin{equation*}
\begin{split}
&D_tu_R(t,x) + \A^tu_R(t,x) - c(t,x)u_R(t,x) = 0, \ (t,x) \in [0,T)\times B(R),\\
&u_R(t,x) = 0, \quad t\in [0,T],\ x\in \partial B(R),\\
&u_R(T,x) = \eta(x) h(x), \quad x\in B(R),
\end{split}
\end{equation*}
where $\eta(x) = \eta(|x|)\in C^\infty_0(\er^d)$ such that $\eta(x) = 1$ for $x\in B(\frac12R)$ and $\eta(x) =0$ outside $B(R)$.

By Assumption \ref{as_161} and Theorem \ref{th_143} $u_R(t,\cdot)$ is in $H^{4+\beta}(B(R))$ for any $\beta \in (0,1)$ and $t\in (0, T)$. Then, we define for $a >0$
\begin{align*}
v_R(t,x) = &|u_R(t,x)|^2 + a(T-t)\eta^2\sum_{i=1}^d \big|D^iu_R(t,x)\big|^2 \\
&+ a^2(T-t)^2\eta^4\sum_{i,j=1}^d \big|D^{ij}u_R(t,x)\big|^2 \\
&+ a^3(T-t)^3\eta^6\sum_{i,j,k=1}^d \big|D^{ijk}u_R(t,x)\big|^2.
\end{align*}

Straightforward computations show that the function $v_R$ solves the
Cauchy problem
\begin{equation*}
\begin{split}
&D_tv_R(t,x) + \A^tv_R(t,x) - c(t,x)v_R(t,x) = f_R(t,x), \, (t,x) \in [0,T)\times B(R),\\
&v_R(t,x) = 0, \quad t\in [0,T],\ x\in \partial B(R),\\
&v_R(T,x) = \eta^2(x)h^2(x), \quad x\in B(R),
\end{split}
\end{equation*}
where
\[
f_R(t,x) = \sum_{i=1}^9 f_{i,R}(t,x).
\]
For $i = 1,\dots,8$, functions $f_{i,R}(t,x)$ are essentially equal $-g_{i,n}$, where $g_{i,n}$ are functions in the proof of Theorem 2.4 in \cite{L2011}. The only difference is that terms $(t-s)^k$ in the original formulas for $g_{i,n}$ are replaced by $(T-t)^k$ in $f_{i,R}$. Since Hypotheses 2.1 of \cite{L2011} concerning estimates of the coefficients of equation \eqref{eq_161} are fulfilled by our Assumption \ref{as_161}, then also the estimate (2.23) of \cite{L2011} is valid (time regularity of the coefficients is not important in these estimates, hence our weaker assumptions do not play a role) and
\[
f'_R = \sum_{i=1}^8 f_{i,R}(t,x)\ge 0.
\]
For $f_{9,R}(t,x)$ tedious computations give
\[
f_{9,R} = c(t,x)v_R + g_R
\]
where
\begin{align*}
g_R = &2a(T-t)\eta^2 u_R\sum_{i=1}^d D^icD^iu_R\\
&+ 4a^2(T-t)^2 \eta^4 \sum_{i,j=1}^d D^{ij}u_R D^iu_RD^jc \\
&+ 2a^2(T-t)^2 \eta^4 u_R \sum_{i,j=1}^d D^{ij}u_R D^{ij}c\\
&+ 2 a^3 (T-t)^3 \eta^6 u_R \sum_{i,j, k=1}^d D^{ijk}u_R D^{ijk}c\\
&+ 6a^3(T-t)^3\eta^6 \sum_{i,j,k=1}^d D^icD^{jk}u_R D^{ijk}u_R \\
&+ 6a^3(T-t)^3\eta^6\sum_{i,j,k=1}^d D^{jk}cD^{i}u_R D^{ijk}u_R.
\end{align*}

Applying Cauchy's inequality, we get
\begin{align*}
|g_R| \le C\, \bigl(1+c(t,x)\bigr) &\biggl( 2aT \eta^2 \sum_{i=1}^d \Bigl(\frac{\varepsilon}{2}(D^iu_R)^2 + \frac{1}{2\varepsilon}u^2_R\Bigr)\biggr.\\
&+ 4a^2T^2 \eta^4\sum_{i,j=1}^d \Bigl(\frac{\varepsilon}{2}(D^{ij}u_R)^2 + \frac{1}{2\varepsilon}(D^iu_R)^2\Bigr)\\
&+ 2a^2T^2\eta^4 \sum_{i,j=1}^d \Bigl(\frac{\varepsilon}{2}(D^{ij}u_R)^2 + \frac{1}{2\varepsilon}( u_R)^2\Bigr)\\
&+ 2a^3T^3 \eta^6\sum_{i,j,k=1}^d \Bigl(\frac{\varepsilon}{2}(D^{ijk}u_R)^2 + \frac{1}{2\varepsilon}(u_R)^2\Bigr)\\
&+ 6a^3T^3 \eta^6\sum_{i,j,k=1}^d \Bigl(\frac{\varepsilon}{2}(D^{ijk}u_R)^2 + \frac{1}{2\varepsilon}(D^{ij}u_R)^2\Bigr)\\
&\biggl.+ 6a^3T^3 \eta^6 \sum_{i,j,k=1}^d \Bigl(\frac{\varepsilon}{2}(D^{ijk}u_R)^2 + \frac{1}{2\varepsilon}(D^{i}u_R)^2\Bigr)\biggr).
\end{align*}
Taking $\varepsilon$ sufficiently small and $a\le a_0$, we have
\begin{align*}
C \varepsilon \Bigl(a_0 T \eta^2\sum_{i =1}^d (D^iu_R)^2 &+ 3 a_0^2T^2 \eta^4\sum_{i,j=1}^d (D^{ij}u_R)^2 \Bigr.\\
\Bigl.&+ 7 a_0^3T^3 \eta^6\sum_{i,j,k=1}^d (D^{ijk}u_R)^2 \Bigr) \le v_R/4.
\end{align*}
Then we can select $a$ so small that
\begin{align*}
\frac{CaT\eta^2}{\varepsilon}\Bigl(\!(d+ a_0Td^2\eta^2 &+ a_0^2 T^2 d^3\eta^4) u_R^2 + (2a_0T d\eta^2+ 3a_0^2T^2 d^2\eta^4) \sum_{i =1}^d (D^iu_R)^2\Bigr. \\
\Bigl.&+ 3 a_0^2T^2 d \eta^4 \sum_{i,j=1}^d (D^{ij}u_R)^2\Bigr)\le v_R/4.
\end{align*}

Collecting these estimates, we get
\[
|g_R(t,x)| \le \bigl(1+c(t,x)\bigr) v_R(t,x)/2
\]
and
\[
f_{9,R}(t,x) \ge c(t,x)v_R - \bigl(1+ c(t,x)\bigr) v_R /2 = c(t,x)v_R /2 - v_R/2.
\]

Since
\[
c(t,x)v_R(t,x) + f_R(t,x) \ge (3c_0/2 - 1/2)v_R(t,x) = c_1 v_R(t,x),
\]
then applying Theorem \ref{th_162}, we obtain
\[
\sup_{x\in\er^d} |v_R(t,x)| \le \e^{-c_1(T-t)} \sup_{x\in\er^d} |\eta^2(x)h^2(x)| \le C\,\sup_{x\in\er^d} | h^2(x)|, \quad t\in[0,T].
\]

As the constant in the above estimate is independent of $R$, we can pass to the limit $R\to\infty$ obtaining the estimate of the theorem. Let us observe that under Assumption \ref{as_161} $u(t,x)$, a solution of \eqref{eq_161}, has a continuous third spatial derivative and $D^3_xu_R$ converges to $D^3_xu$ locally uniformly on $[0,T)\times \er^d$. This convergence follows from the observation that fourth derivatives of $u_R(t,x)$ are uniformly bounded on any compact set $\mathcal{B}\subset [0,T)\times \er^d$. Then $D^3_xu_R$ is a sequence of bounded and equicontinuous functions and the Arzel\'a-Ascoli theorem gives the convergence.
This observation is essential for the interpolation argument used in \cite{L2011} to extend the estimate to $\beta_1 \neq 0$ and/or $\beta_2 \neq 0$.
\qed
\medskip

The following theorem is a part of Theorem 3.8 of \cite{L2011} where the author assumes the spatial H\"older continuity and only the continuity in time of the coefficients of \eqref{eq_142}.

\bt\label{th_164}
Let the conditions of Assumption \ref{as_161} hold and $u(t,x)$ be a solution of the Cauchy problem \eqref{eq_142} with $h=0$ and $f(t,\cdot) \in H^{\beta} (\er^d)$, $\beta\in (0,1)$. Then
\[
\| u(t, \cdot)\|_{H^{2+\beta}(\er^d)} \le C \sup_{t\in[0,T]} \| f(t, \cdot)\|_{H^{\beta}(\er^d)}, \quad t\in [0,T).
\]
\et

The following corollary summarizes the results of the two previous theorems.
\begin{cor}\label{cor_161}
Let Assumption \ref{as_161} hold and $f(t,\cdot) \in H^{\beta} (\er^d)$ for $t\in[0,T]$, $\beta\in (0,1)$. Then for a solution $u(t,x)$ of the Cauchy problem \eqref{eq_142} we have the estimate
\begin{align*}
\| u(t, \cdot)\|_{BC(\er^d)} &+ (T-t)^{1+\beta/2}\| u(t, \cdot)\|_{H^{2 +\beta}(\er^d)} \\
&\le C \Bigl(\|h\|_{BC(\er^d)} + \sup_{t\in[0,T]} \| f(t, \cdot)\|_{H^{\beta}(\er^d)}\Bigr), \quad t\in [0,T).
\end{align*}

\end{cor}

From Theorem 3.8 of \cite{L2011}, we obtain the optimal Schauder's estimate.
\bt\label{th_165}
Let the conditions of Assumption \ref{as_161} hold and $u(t,x)$ be a solution of the Cauchy problem \eqref{eq_142} with $h\in H^{2+\beta} (\er^d)$ and $f(t,\cdot) \in H^{\beta} (\er^d)$ for $t\in [0,T]$, $\beta\in (0,1)$. Then
\[
\sup_{t\in[0,T]} \| u(t, \cdot)\|_{H^{2 +\beta}(\er^d)} \le C\Bigl(\|h\|_{H^{2 +\beta}(\er^d)} + \sup_{t\in[0,T]} \| f(t, \cdot)\|_{H^{\beta}(\er^d)}\Bigr).
\]
\et

We now prove the main result of this section: Schauder's estimates in the weighted H\"older spaces with polynomial weights for the Kolmogorov equation \eqref{eq_142}.

\bt\label{th_166}
Let Assumption \ref{as_161} hold, $h\in BC_P(\er^d)$, and $f(t,\cdot) \in H_P^{\beta} (\er^d)$ for $t\in[0,T]$, $\beta\in (0,1)$, where $P(x) = 1 +|x|^{2q}$, $q\ge 0$. In addition, we assume that the coefficients $a_{ij}(t,x)$ are four times continuously differentiable with respect to $x$ and $D^\alpha_xa_{ij}(t,x)\in H^\beta_{\mathit{loc}}(\er^d)$ for $|\alpha| = 4$. Then for a solution $u(t,x)$ of the Cauchy problem \eqref{eq_142} we have the estimate
\begin{align*}
\| u(t, \cdot)\|_{BC_P(\er^d)} &+ (T-t)^{1+\beta/2}\| u(t, \cdot)\|_{H_P^{2 +\beta}(\er^d)} \\
&\le C \Bigl(\|h\|_{BC_P(\er^d)} + \sup_{t\in[0,T]} \| f(t, \cdot)\|_{H_P^{\beta}(\er^d)}\Bigr), \quad t\in [0,T).
\end{align*}
If in addition, $h\in H_P^{2+\beta} (\er^d)$, then we have the estimate
\[
\sup_{t\in[0,T]} \| u(t, \cdot)\|_{H_P^{2 +\beta}(\er^d)} \le C\Bigl(\|h\|_{H_P^{2 +\beta}(\er^d)} + \sup_{t\in[0,T]} \| f(t, \cdot)\|_{H_P^{\beta}(\er^d)}\Bigr).
\]
\et
\proof
To prove the theorem, we transform the Cauchy problem \eqref{eq_142} with un\-boun\-ded $h(x)$ and $f(t,x)$ into a problem with bounded $h(x)$ and $f(t,x)$ where Schau\-der's estimates are given in Corollary \ref{cor_161} and Theorem \ref{th_165}.

Assume that $u(t,x)$ is a solution of \eqref{eq_142} with unbounded $h(x)$ and $f(t,x)$. Let $u(t,x) = v(t,x)P(x)$. Then $v(t,x)$ fulfills the equation
\begin{equation}\label{eq_163}
\begin{split}
&D_tv(t,x) + \tilde\A^tv(t,x) - \tilde c(t,x)v(t,x) + \tilde f(t,x) =0, \ (t,x) \in [0,T)\times\er^d,\\
&v(T,x) = \tilde h(x), \quad x\in \er^d,
\end{split}
\end{equation}
where
\[
\tilde\A^tv(t,x) = \sum_{i,j =1}^d a_{ij}(t,\bbx)D^{ij}_x v(t,x) + \sum_{i=1}^d \tilde b_i(t,\bbx) D^i_x v(t,x),
\]
and
\begin{align*}
&\tilde b_i(t,x) = b_i(t,x) + \sum_{j=1}^d D^j_x a_{ij}(t,x)P^{-1}(x)D^i_xP(x),\\
&\tilde c(t,x) = c(t,x) + \sum_{i,j =1}^d a_{ij}(t,\bbx)P^{-1}(x) D^{ij}_x P(x) + \sum_{i=1}^d b_i(t,\bbx)P^{-1}(x) D^i_x P(x),\\
& \tilde f(t,x) = P^{-1}(x)f(t,x),\\
& \tilde h(x) = P^{-1}(x) h(x).
\end{align*}

It is obvious that $\tilde h\in BC(\er^d)$ and $\tilde f(t, \cdot) \in H^\beta(\er^d)$. To show that $\tilde c(t,x)$ fulfills the conditions of Lemma \ref{lem_161} let us observe that
\[
| a_{ij}(t,\bbx)P^{-1}(x) D^{ij}_x P(x)| + |b_i(t,\bbx)P^{-1}(x) D^i_x P(x)| \le C,
\]
which follows from the estimate
\[
|P^{-1}(x) D_x^\alpha P(x)|\le C (1+|x|^2)^{-|\alpha|/2}
\]
and the growth conditions for $a_{ij}$ and $b_i$.

Then
\begin{align*}
\tilde c(t,x) \ge c_0 &- \sup_{x\in \er^d} \Big|\sum_{i,j =1}^d a_{ij}(t,\bbx)P^{-1}(x) D^{ij}_x P(x)\Big| \\
&- \sup_{x\in\er^d} \Big|\sum_{i=1}^d b_i(t,\bbx)P^{-1}(x) D^i_x P(x)\Big| \ge \tilde c_0.
\end{align*}

Since from Assumption \ref{as_161} we have
\[
|D^\alpha_x a_{ij}(t,x)| \le C(1+|x|), \quad |D^\alpha_x b_{i}(t,x)| \le C,
\]
for $1 \le |\alpha| \le 3$, and
\[
\big|D^{\alpha_2}_x \bigl(P^{-1}(x) D^{\alpha_1}_x P(x)\bigr)\big| \le C(1+|x|^2)^{-|\alpha_1|/2-|\alpha_2|/2},
\]
then for $|\alpha| \le 3$
\begin{align*}
\Big|D^\alpha_x \Bigl( \sum_{i,j =1}^d a_{ij}(t,\bbx)P^{-1}(x) D^{ij}_x P(x) &+ \sum_{i=1}^d b_i(t,\bbx)P^{-1}(x) D^i_x P(x)\Bigr)\Big|
\le C.
\end{align*}
Therefore, $\tilde c(t,x)$ fulfills condition 3 of Assumption \ref{as_161}.

From the above estimates, we also obtain
\[
\big|D^j_x a_{ij}(t,x)P^{-1}(x)D^i_xP(x)\big| \le C,
\]
and for $|\alpha| \le 3$
\[
\Big|D^\alpha_x \Bigl(\sum_{j=1}^d D^j_x a_{ij}(t,x)P^{-1}(x)D^i_xP(x)\Bigr)\Big|\le C.
\]
The additional assumption that $a_{ij}(t,x)$ is four times differentiable with respect to $x$ and $D^\alpha_xa_{ij}(t,x)\in H^\beta(\er^d)$ for $|\alpha| = 4$ guaranties that $D^\alpha_x\tilde b_i(t,x) \in H^\beta(\er^d)$ for $|\alpha| =3$. Under this assumption condition 1 of Assumption \ref{as_161} is fulfilled by $\tilde b(t,x)$.

Thus, problem \eqref{eq_163} fulfills the conditions of Corollary \ref{cor_161} and Theorem \ref{th_165} and the relevant Schauder's estimates for $v(t,x)$ hold. This proves the estimates of the theorem for $u(t,x) = v(t,x)P(x)$.
\qed

\end{document}